\documentclass{article}
\usepackage{amsfonts}
\usepackage{amssymb}
\usepackage{amsmath}
\usepackage{amsthm}
\usepackage{euscript}
\usepackage[dvips]{graphics}

\newtheorem{theorem}{Theorem}[section]
\newtheorem{proposition}[theorem]{Proposition}

\newcommand{\bu}{ {\rm bulk}}

   \title{Shuffling algorithm for boxed plane partitions}
   \author{Alexei Borodin\thanks{California Institute of Technology and Institute for Information Transmission Problems,
   Moscow}  \ and
   Vadim Gorin\thanks{Moscow State University} }
\date{}

\begin{document}
\maketitle

\begin{abstract}
We introduce discrete time Markov chains that preserve uniform
measures on boxed plane partitions. Elementary Markov steps change
the size of the box from $a\times b \times c$ to $(a-1)\times
(b+1)\times c$ or $(a+1)\times (b-1)\times c$. Algorithmic
realization of each step involves $O((a+b)c)$ operations. One
application is an efficient perfect random sampling algorithm for
uniformly distributed boxed plane partitions.

Trajectories of our Markov chains can be viewed as random point
configurations in the three-dimensional lattice. We compute the bulk
limits of the correlation functions of the resulting random point
process on suitable two-dimensional sections. The limiting
correlation functions define a two-dimensional determinantal point
processes with certain Gibbs properties.
\end{abstract}

\section*{Intoduction}

For any integers $a,b,c\ge 1$ consider a hexagon with sides
$a,b,c,a,b,c$ drawn on a regular triangular lattice. Denote by
$\Omega_{a\times b\times c}$ the set of all tilings of this hexagon
by rhombi obtained by gluing two of the neighboring elementary
triangles together (such rhombi are called {\it lozenges\/}).
Equivalently, $\Omega_{a\times b\times c}$ is the set of all dimers
on the part of the dual hexagonal lattice cut out by our
$(a,b,c)$-hexagon. An element of $\Omega_{4\times 5\times 5}$ is
shown on Figure 1.

Elements of $\Omega_{a\times b\times c}$ have a number of different
interpretations, see e.g. Section 1 below. In particular, they can
be viewed as plane partitions or stepped surfaces inside a
three-dimensional box of size $a\times b\times c$. A uniformly
distributed element of $\Omega_{a\times b\times c}$ then provides a
basic model of a {\it random surface\/}. This model and its
generalizations have been thoroughly studied, see e.g.
\cite{CohnLarsenPropp}, \cite{CohnKenyonPropp}, \cite{DMB},
\cite{Des}, \cite{J1} \cite{J}, \cite{JN}, \cite{KenyonHeight},
\cite{KenyonOkounkov}, \cite{Kr}, \cite{Propp1}, \cite{Propp2},
 \cite{W1}, \cite{W2}.

\begin{center}
 {\scalebox{0.8}{\includegraphics{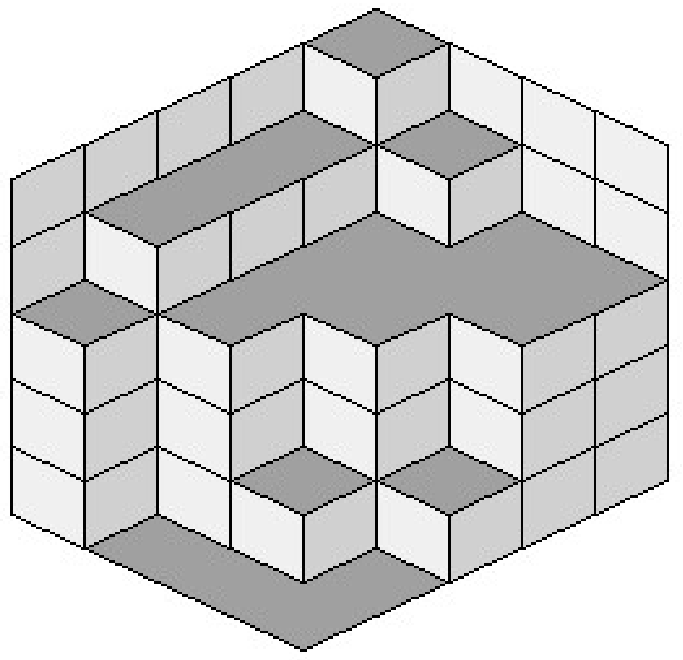}}}

 Figure 1. Tiling, stepped surface or boxed plane partition. Separating
 lines between horizontal lozenges are removed.
\end{center}

The main goal of this paper is to introduce and study certain
discrete time Markov chains on boxed plane partitions that preserve
the uniform measures.

Denote by $\mu_{a\times b\times c}$ the uniform probability measure
on $\Omega_{a\times b\times c}$. We construct two families of
stochastic matrices
$$
\gathered P_{a\times b\times c}^+:\Omega_{a\times b\times c} \times
\Omega_{(a-1)\times (b+1)\times c}\to [0,1], \\ P^-_{a\times b\times
c}:\Omega_{a\times b\times c} \times \Omega_{(a+1)\times (b-1)\times
c}\to [0,1],
\endgathered
$$
such that
$$
\sum_{\omega\in\Omega_{a\times b\times c}}\mu_{a\times b\times
c}(\omega)P_{a\times b\times c}^\pm(\omega,\omega')=\mu_{(a\mp
1)\times (b\pm1)\times c}(\omega')
$$
for all $a,b,c$ such that the participating sets $\Omega$ are
nonempty.

Although it is a little awkward to write matrices $P^\pm$ in one
formula, application of Markov operators corresponding to these
matrices is fairly easy to describe algorithmically. The exact
algorithm can be found in Section \ref{Section_Algorithm}. Roughly
speaking, the algorithm does the following: Given
$\omega\in\Omega_{a\times b\times c}$, in order to construct
$\omega'$ distributed according to $P^\pm_{a\times b\times
c}(\omega, \omega')$ it needs to consider all horizontal lozenges of
$\omega$ sequentially from left to right. For each such lozenge the
algorithm decides on its new position using a simple one-dimensional
probability distribution. (Here we ignored occasional
appearance/diappearance of horizontal lozenges on top and at the
bottom, that sometimes take place.) In a way, this means that
$P^\pm$ decompose into products of one-dimensional Markov steps.

The algorithm has some similarity to the shuffling algorithm for
domino tilings of the Aztec diamonds introduced in \cite{EKLP}.
Indeed, that algorithm also maps uniform measures to uniform
measures (actually, it works for a one-dimensional family of
measures that includes the uniform one), and it also decomposes into
simple (Bernoulli) Markov steps. For that reason we call our
algorithm the {\it shuffling algorithm for boxed plane
partitions\/}.

We further consider Markov chains obtained by successive application
of arbitrary sequences of matrices $P^+$ and $P^-$. The initial
condition and the one-time distributions are all uniform measures on
the appropriate spaces $\Omega$. One example is the application of
the sequence of $b$ matrices $P^+$ to the unique probability measure
on singleton $\Omega_{(a+b)\times 0\times c}$. This gives a perfect
sampling algorithm for $\mu_{a\times b\times c}$ that takes
$O((a+b)bc)$ one-dimensional steps. When $a$, $b$, and $c$ are
comparable, the algorithm is roughly as efficient as that from
\cite{Kr}, and it is more efficient than other known algorithms, cf.
\cite{BM}, \cite{LRS}, \cite{Propp1}, \cite{Propp2}, \cite{W1},
\cite{W2}. Another example is an alternating sequence
$(P^+P^-)(P^+P^-)\cdots$ that provides an equilibrium dynamics on
$\Omega_{a\times b\times c}$ with $\mu_{a\times b\times c}$ as the
equilibrium measure.

These Markov chains can be viewed as two-dimensional random
growth/decay models. One important feature of these models is the
fact that on certain two-dimensional sections of the
three-dimensional space-time, their suitably defined correlation
functions are computable in a closed determinantal form. For the
static correlation functions (those of the measures $\mu_{a\times
b\times c}$) this fact is well known, see \cite{Kasteleyn},
\cite{KenyonLocal}, \cite{J}, \cite{Gor}.

We then focus on the {\it bulk asymptotics\/} of the correlation
functions. Namely, when $a,b,c\to\infty$ with $a:b:c$ fixed, we look
at the three-dimensional lattice process near a fixed point of the
global limit shape of our random surface. On the same
two-dimensional sections we compute the limiting correlation
functions that also have the determinantal form. In fact, for every
such section, they define a two-dimensional determinantal random
point process with a certain Gibbs property, see \cite{BS} and
references therein. Our results extend similar results for the
static case obtained in \cite{Gor}.

The construction of $P^\pm$ is an application of the general
algebraic formalism developed in \cite{BF}. The continuous time
Markov chain considered in \cite{BF} in detail can be viewed as the
degeneration of $P^\pm$ near a corner of the hexagon as $a,b,c$
become large, and either $a$ and $b$ is substantially larger than
the other two. It is worth noting that the shuffling algorithm for
domino tilings of the Aztec diamonds also fits into the formalism of
\cite{BF}, see Section 2.6 of \cite{BF} and \cite{Nordenstam}.

The proof of the bulk asymptotics involves spectral decomposition of
the matrices $P^\pm$ in terms of Hahn classical orthogonal
polynomials. However, the limiting argument does not require the
asymptotics of Hahn polynomials themselves --- it only requires the
much simpler asymptotics of the difference operators related to Hahn
polynomials. This approach to bulk asymptotics of determinantal
point processes is due to G.~Olshanski; it was first used in
\cite{BO2} for Charlier and Krawtchouk orthogonal polynomials, and
in \cite{Gor} it was further developed in the more complex case of
Hahn polynomials.

The paper is organized as follows. In Section 1 we discuss different
combinatorial definitions of our model and introduce notations. In
Section 2, we introduce four stochastic matrices on one-dimensional
point configurations and prove certain commutativity relations
between them. In Section 3, we use these four matrices to define
$P^\pm$; the construction is based on an idea from \cite{DF}.
Section 4 contains the algorithmic description of $P^\pm$ and images
obtained from their computer realizations. In Section 5 we compute
the correlation functions on suitable two-dimensional sections of
the space-time, and in Section 6 we obtain the bulk limits.

The first named author (A. B.) was partially supported by the NSF
grant DMS-0707163. The second named author (V. G.) was partially
supported by RFBR grant 07-01-91209, the Moebius Contest Foundation
for Young Scientists and Leonhard Euler's Fund of Russian
Mathematics Support.

\section{Basic model}
The main object of our study has many different combinatorial
interpretations. In this section we discuss some of them.

Consider a tiling of an equi-angular hexagon of side lengths $a$,
$b$, $c$, $a$, $b$, $c$ by rhombi with angles $\pi/3$ and $2\pi/3$
and side lengths $1$. Such rhombi are called lozenges.

Lozenge tilings of a hexagon can be identified with 3d Young
diagrams (equivalently, boxed plane partitions) or with stepped
surfaces. The bijection is best described pictorially. Examine
Figure 1, where a tiling of the $(4\times 5\times 5)$ hexagon is
shown, and view this picture as a 3d shape.

Given a tiling we construct a family of non-intersecting paths on
the surface of the corresponding 3d Young diagram as in \cite{LRS},
\cite{J}. Figure 2 provides an example.

\begin{center}
 {\scalebox{0.8}{\includegraphics{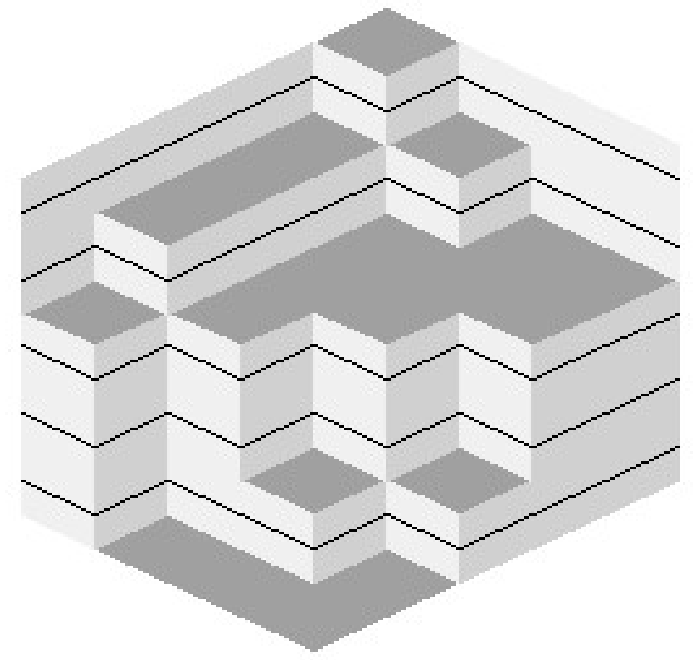}}}

 Figure 2. Non-intersecting paths on the surface of 3d Young diagram.
\end{center}

We view this family as a family of paths on the plane. It is
convenient for us to do one more modification. We replace the
downgoing segments of paths by horizontal ones and we replace
upgoing segments by segments of slope 1. Consequently, our family is
interpreted as a family of non-intersecting paths on $\mathbb Z^2$.
Figure 3 shows the family corresponding to the tiling on Figures 1
and 2.

\begin{center}
 {\scalebox{0.8}{\includegraphics{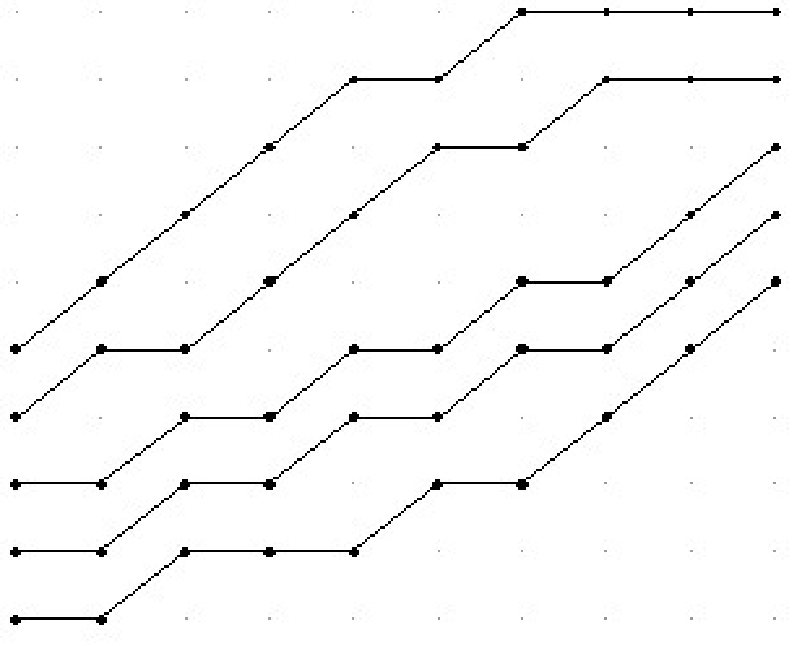}}}

 Figure 3. Non-intersecting paths on $\mathbb Z^2$.
\end{center}

Below we are going to use the last interpretation.

Let us introduce some notations. Denote by $N$ the number of paths
(in our example $N=5$). Introduce coordinates on $\mathbb Z^2$ so
that the first path starts at the point $(0,0)$ and ends at the
point $(T,S)$. The second path starts at the point $(0,1)$ and ends
at the point $(T,S+1)$, and so on. Finally, the $N$th path starts at
the point $(0,N-1)$ and ends at the point $(T,S+N-1)$. In our
example $T=9$ and $S=5$.

Denote by $\Omega(N,T,S)$ the set of families of $N$
non-intersecting paths made of segments of slopes 0,1, starting at
$(0,0),\dots,(0,N-1)$ and ending at $(T,S),\dots,(T,S+N-1)$. Note
that $\Omega(N,T,S)\neq \emptyset$, is equivalent to $0\le S\le T$.

 Denote by $\mu(N,T,S)$ the uniform
measure on $\Omega(N,T,S)$.

Note that $\Omega(N,T,S)$ was called $\Omega_{a\times b\times c}$ in
Introduction, while  measure $\mu(N,T,S)$ was called $\mu_{a\times
b\times c}$. Here $a=T-S$, $b=S$, $c=N$.

Set
$$
 \mathfrak X_{N,T}^{S,t}=\{x\in\mathbb Z: \max(0,t+S-T)  \le x\le \min(t+N-1,S+N-1)\}
$$
and
$$
 \mathcal X_{N,T}^{S,t}=\{ (x_1,x_2,\dots,x_N)\in (\mathfrak
 X_{N,T}^{S,t})^N: x_1<x_2<\dots<x_N\}.
$$
$\mathfrak X_{N,T}^{S,t}$ is the section of our hexagon by the
vertical line with coordinate $t$, and $\mathcal X_{N,T}^{S,t}$ is
the set of all $N$--tuples  in this section.

For any  $\mathcal T\in\Omega(N,T,S)$ denote by
$\tau_1,\tau_2,\dots,\tau_N$ the corresponding paths, numbering
starts from the bottom one. Thus,
$$
\tau_i= \{i-1=\tau_i(0),\tau_i(1),\dots,\tau_i(T-1),\tau_i(T)=S+i-1
\},
$$
so that $\tau_i(t+1)-\tau_i(t)\in\{0,1\}$. Note that
$$\tau_i(t)\in \mathfrak X_{N,T}^{S,t},\quad
\left( \tau_1(t),\dots,\tau_N(t) \right) \in \mathcal
X_{N,T}^{S,t}.$$

Consequently, any family of paths $\mathcal T$ can be identified
with a sequence
$$\{X(1),\dots,X(T)\},\quad X(t)\in \mathcal X_{N,T}^{S,t},$$
where
$$
 X(t)=\left(\tau_1(t),\dots, \tau_N(t) \right).
$$

In fact, $X(t)$ is a Markov chain with time $t$, see \cite{J1},
\cite{J}, \cite{JN}, \cite{Gor}. Its transition probabilities are
given in Proposition \ref{Proposition_tr_prob} below.

Through the Sections
\ref{Section_Four_families}-\ref{Section_General_process} the
parameters $N$ and $T$ remain fixed and we omit them in different
notations. We write $\mathfrak X^{S,t}$ instead of $\mathfrak
X_{N,T}^{S,t}$, $\mathcal X^{S,t}$ instead of $\mathcal
X_{N,T}^{S,t}$, and so on.

In the present paper we introduce a discrete time Markov chain
$M(r)$, where $r$ is a time variable. $M(r)$ takes values in
$\Omega(N,T,r)$ and one-dimensional distributions of $M(r)$ coincide
with $\mu(N,T,r)$.

\section{Four families of stochastic matrices}
\label{Section_Four_families}
\subsection{Properties of sections $X(t)$}

Denote by $\rho_{S,t}$ the projection of $\mu(N,T,S)$ to $\mathcal
X^{S,t}$, i.e.
$$
\rho_{S,t}(Y)={\rm Prob}\{X(t)=Y\},\quad
Y=(y_1<\dots<y_N)\in\mathcal X^{S,t}.
$$

The following two propositions were proved in \cite[Theorem 4.1]{J1}
and \cite[Lemma 4]{Gor}. Below we use the Pochhammer symbol
$(a)_k=\Gamma(a+k)/\Gamma(a)=a(a+1)\cdots (a+k-1)$.

\begin{proposition}
\label{Proposition_measure} We have
$$\rho_{S,t}(Y)=Z_{S,t}
\prod\limits_{1\le i<j\le N}(y_i-y_j)^2\prod\limits_{i=1}^N
w^{S,t}(y_i),$$ where
$$
w^{S,t}(x)=\frac{1}{x!(t+N-1-x)!(S+N-1-x)!(T-t-S+x)!}
$$ and
$$
Z_{S,t}=\prod\limits_{i=1}^N\frac{(t+1)_{i-1}(T-t+1)_{i-1}(S-i+N)!(T-S+i-1)!}
    {(T+1)_{i-1} (i-1)!}
  \cdot \left(\frac{t!(T-t)!}{T!}\right)^N.$$
\end{proposition}

\begin{proposition}
\label{Proposition_tr_prob} For $t=0,1,\dots,T-1,$
\begin{multline*}
{\rm Prob}\{X(t+1)=(y_1,\dots ,y_N)\mid X(t)=(x_1,\dots ,x_N)\}\\=
 \frac{\prod\limits_{i<j}(y_j-y_i)\cdot
   \prod\limits_{i:\,y_i=x_i+1}(N+S-x_i-1)\cdot\prod\limits_{i:\,y_i=x_i}(x_i+T-t-S)}
 {(T-t)_{N}\cdot\prod\limits_{i<j}(x_j-x_i)},
\end{multline*}
if  $y_i-x_i\in\{0,1\}$ for every $i$, and the conditional
probability is equal to zero otherwise.
\end{proposition}

We also need ``co-transition probabilities'' of the Markov chain
$X(t)$.
\begin{proposition}
\label{Proposition_co-tr} For $t=1,2,\dots,T$
\begin{multline*}
{\rm Prob}\{X(t-1)=(y_1,\dots ,y_N)\mid X(t)=(x_1,\dots ,x_N)\}\\=
 \frac{\prod\limits_{i<j}(y_j-y_i)\cdot
   \prod\limits_{i:\,y_i=x_i-1}x_i\cdot\prod\limits_{i:\,y_i=x_i}(t+N-1-x_i)}
 {(t)_{N}\cdot\prod\limits_{i<j}(x_j-x_i)},
\end{multline*}
if  $y_i-x_i\in\{-1,0\}$ for every $i$, and the conditional
probability is equal to zero otherwise.
\end{proposition}
\begin{proof}
Straightforward computation using
\begin{multline*}
 {\rm Prob}\{X(t-1)=(y_1,\dots ,y_N)\mid
X(t)=(x_1,\dots ,x_N)\}\\= {\rm Prob}\{X(t)=(x_1,\dots ,x_N)\mid
X(t-1)=(y_1,\dots ,y_N)\}
\cdot\frac{\rho_{S,t-1}(y_1,\dots,y_N)}{\rho_{S,t}(x_1,\dots,x_N)}
\end{multline*}
\end{proof}

\subsection{Stochastic matrices}

We need four families of stochastic matrices $P^{S,t}_{t+}$,
$P^{S,t}_{t-}$, $P^{S,t}_{S+}$, $P^{S,t}_{S-}$.

$P^{S,t}_{t+}(X,Y)$ is an $|\mathcal X^{S,t}|\times|\mathcal
X^{S,t+1}|$ matrix, $X=(x_1<\dots<x_N)\in\mathcal X^{S,t}$,
$Y=(y_1<\dots<y_n)\in\mathcal X^{S,t+1}$;

$$
P^{S,t}_{t+}(X,Y)= \frac{\prod\limits_{i<j}(y_j-y_i)
   \prod\limits_{i:\,y_i=x_i+1}(N+S-x_i-1)\prod\limits_{i:\,y_i=x_i}(x_i+T-t-S)}
 {(T-t)_{N}\cdot\prod\limits_{i<j}(x_j-x_i)},
$$
if  $y_i-x_i\in\{0,1\}$ for every $i$, and $P^{S,t}_{t+}(X,Y)=0$
otherwise.

$P^{S,t}_{S+}(X,Y)$ is an $|\mathcal X^{S,t}|\times|\mathcal
X^{S+1,t}|$ matrix, $X=(x_1<\dots<x_N)\in\mathcal X^{S,t}$,
$Y=(y_1<\dots<y_n)\in\mathcal X^{S+1,t}$;

$$
P^{S,t}_{S+}(X,Y)= \frac{\prod\limits_{i<j}(y_j-y_i)
   \prod\limits_{i:\,y_i=x_i+1}(N+t-x_i-1)\prod\limits_{i:\,y_i=x_i}(x_i+T-t-S)}
 {(T-S)_{N}\cdot\prod\limits_{i<j}(x_j-x_i)},
$$
if  $y_i-x_i\in\{0,1\}$ for every $i$, and $P^{S,t}_{S+}(X,Y)=0$
otherwise.

$P^{S,t}_{t-}(X,Y)$ is an $|\mathcal X^{S,t}|\times|\mathcal
X^{S,t-1}|$ matrix, $X=(x_1<\dots<x_N)\in\mathcal X^{S,t}$,
$Y=(y_1<\dots<y_n)\in\mathcal X^{S,t-1}$;

$$
P^{S,t}_{t-}(X,Y)= \frac{\prod\limits_{i<j}(y_j-y_i)
   \prod\limits_{i:\,y_i=x_i-1}x_i\prod\limits_{i:\,y_i=x_i}(t+N-1-x_i)}
 {(t)_{N}\cdot\prod\limits_{i<j}(x_j-x_i)},
$$
if  $y_i-x_i\in\{-1,0\}$ for every $i$, and $P^{S,t}_{t-}(X,Y)=0$
otherwise.

$P^{S,t}_{S-}(X,Y)$ is an $|\mathcal X^{S,t}|\times|\mathcal
X^{S-1,t}|$ matrix, $X=(x_1<\dots<x_N)\in\mathcal X^{S,t}$,
$Y=(y_1<\dots<y_n)\in\mathcal X^{S-1,t}$;

$$
 P^{S,t}_{S-}(X,Y)= \frac{\prod\limits_{i<j}(y_j-y_i)
   \prod\limits_{i:\,y_i=x_i-1}x_i\prod\limits_{i:\,y_i=x_i}(S+N-1-x_i)}
 {(S)_{N}\cdot\prod\limits_{i<j}(x_j-x_i)},
$$
if  $y_i-x_i\in\{-1,0\}$ for every $i$, and $P^{S,t}_{S-}(X,Y)=0$
otherwise.

Looking at spaces that parameterize rows and columns of these
matrices one can say that $P^{S,t}_{t+}$ increases $t$,
$P^{S,t}_{t-}$ decreases $t$, while $P^{S,t}_{S+}$ increases $S$ and
$P^{S,t}_{S-}$ decreases $S$.

\begin{theorem}
\label{Th_4stochmatr} All four types of matrices defined above are
stochastic, and they preserve the family of measures $\rho_{S,t}$.
In other words
\begin{equation}
\label{Simple_stM}
 \sum\limits_{Y\in\mathcal X^{S,t\pm 1}} P^{S,t}_{t\pm}(X,Y)=1,\quad
 \sum\limits_{Y\in\mathcal X^{S,S\pm 1}} P^{S,t}_{t\pm}(X,Y)=1,
\end{equation}
$$
 \rho_{S,t\pm 1}(Y)=\sum\limits_{X\in\mathcal X^{S,t}}
 P^{S,t}_{t\pm}(X,Y)\cdot\rho_{S,t}(X),\,
 \rho_{S\pm 1,t}(Y)=\sum\limits_{X\in\mathcal X^{S,t}}
 P^{S,t}_{S\pm}(X,Y)\cdot\rho_{S,t}(X).
$$

\end{theorem}
\begin{proof}
Propositions \ref{Proposition_tr_prob} and \ref{Proposition_co-tr}
imply the claim for $P^{S,t}_{t+}(X,Y)$ and $P^{S,t}_{t-}(X,Y)$.

Now observe that the space $\mathcal X^{S,t}$ is unaffected when we
interchange parameters $t$ and $S$, i.e.
$$
 \mathcal X^{S,t}=\mathcal X^{t,S}.
$$

Moreover, the measures $\rho_{S,t}$ are also invariant under
$S\leftrightarrow t$, i.e.
$$
\rho_{S,t}=\rho_{t,S}.
$$

Finally, note that $P^{S,t}_{t+}(X,Y)$ becomes $P^{S,t}_{S+}(X,Y)$
under $S\leftrightarrow t$ and $P^{S,t}_{t-}(X,Y)$ becomes
$P^{S,t}_{S-}(X,Y)$.

Consequently, applying $S\leftrightarrow t$ to the relations for
$P_{t\pm}$ we obtain needed relations for $P_{S\pm}$.
\end{proof}

\subsection{Determinantal representation}
In this section we write our stochastic matrices in a determinantal
form. This representation is very convenient for various
computations.

First, we introduce $4$ new two-diagonal matrices:
$$
U^{S,t}_{t+}(x,y)=\begin{cases}
              N+S-1-x,&\text{if }y=x+1,\\
              T-t-S+x,&\text{if }y=x,\\
              0,&\text{otherwise,}
             \end{cases}
             \qquad x\in\mathfrak X^{S,t},\quad y\in\mathfrak
             X^{S,t+1};
$$

$$
U^{S,t}_{S+}(x,y)=\begin{cases}
              N+t-1-x,&\text{if }y=x+1,\\
              T-t-S+x,&\text{if }y=x,\\
              0,&\text{otherwise,}
             \end{cases}
             \qquad x\in\mathfrak X^{S,t},\quad y\in\mathfrak
             X^{S+1,t};
$$

$$
U^{S,t}_{t-}(x,y)=\begin{cases}
              x,&\text{if }y=x-1,\\
              t+N-1-x,&\text{if }y=x,\\
              0,&\text{otherwise,}
             \end{cases}
             \qquad  x\in\mathfrak X^{S,t},\quad y\in\mathfrak
             X^{S,t-1};
$$

$$
U^{S,t}_{S-}(x,y)=\begin{cases}
              x,&\text{if }y=x-1,\\
              S+N-1-x,&\text{if }y=x,\\
              0,&\text{otherwise,}
             \end{cases}
             \qquad  x\in\mathfrak X^{S,t},\quad y\in\mathfrak X^{S-1,t}.
$$

It is possible to express stochastic matrices $P^{S,t}_{t\pm}$,
$P^{S,t}_{S\pm}$ through  certain minors of the matrices defined
above.

\begin{proposition} We have
$$P^{S,t}_{t+}(X,Y)= \frac{\prod\limits_{i<j}(y_j-y_i)
   \det[U^{S,t}_{t+}(x_i,y_j)]_{i,j=1,\dots,N}}
 {(T-t)_{N}\cdot\prod\limits_{i<j}(x_j-x_i)}$$

$$P^{S,t}_{S+}(X,Y)= \frac{\prod\limits_{i<j}(y_j-y_i)
   \det[U^{S,t}_{S+}(x_i,y_j)]_{i,j=1,\dots,N}}
 {(T-S)_{N}\cdot\prod\limits_{i<j}(x_j-x_i)}$$

$$P^{S,t}_{t-}(X,Y)= \frac{\prod\limits_{i<j}(y_j-y_i)
   \det[U^{S,t}_{t-}(x_i,y_j)]_{i,j=1,\dots,N}}
 {(t)_{N}\cdot\prod\limits_{i<j}(x_j-x_i)}$$

$$P^{S,t}_{S-}(X,Y)= \frac{\prod\limits_{i<j}(y_j-y_i)
   \det[U^{S,t}_{S-}(x_i,y_j)]_{i,j=1,\dots,N}}
 {(S)_{N}\cdot\prod\limits_{i<j}(x_j-x_i)}$$
\end{proposition}
\begin{proof}
Straightforward computation using the definitions of stochastic
matrices $P^{S,t}_{t\pm}$, $P^{S,t}_{S\pm}$ and matrices
$U^{S,t}_{t\pm}$, $U^{S,t}_{S\pm}$.

Any submatrix of a two-diagonal matrix, which has a nonzero
determinant, is block-diagonal, where each block is either upper or
lower triangular matrix. Thus, any nonzero minor is a product of
suitable matrix elements.
\end{proof}

\subsection{Spectral decomposition of stochastic matrices}
\label{Section_Spectral_decomp} In this section we modify the
determinantal representation of the stochastic matrices and
introduce new representation, which we call \emph{spectral
decomposition}. Spectral decomposition is of crucial importance for
computing correlation functions of the processes that will be
constructed later on, and for finding bulk limits of the processes.
Results of this section will be used in Sections
\ref{Section_General_process} and \ref{Section_Bulk}, while Sections
\ref{Section_M_steps} and \ref{Section_Algorithm} are independent of
these results.

Let us introduce some notations.

Denote by  $H_k^{S,t}(x)$ the Hahn polynomial of the degree $k$
corresponding to the parameters $N$, $T$, $t$, $S$. Domain of
definition of these polynomials coincides with $\mathfrak X^{S,t}$,
and the polynomials are orthogonal with respect to the weight
function $w^{S,t}(x)$ defined in Proposition
\ref{Proposition_measure}. For definition and explicit formulas for
Hahn polynomials see \cite{KS}, more information about the usage of
Hahn polynomials in our case can be found in \cite{Gor}.

Let us denote by $\Psi_k^{S,t}(x)$ the normalized Hahn polynomials
$$\Psi_k^{S,t}(x)=\frac{H_k^{S,t}(x)\sqrt{w^{S,t}(x)}}{\sqrt{(H_k^{S,t},H_k^{S,t})}}.$$
Here $(H_k^{S,t},H_k^{S,t})$ is the squared norm of the polynomial
$H_k^{S,t}$ in $L_2(\mathfrak X^{S,t},w^{S,t}(x))$. The functions
$\Psi_k^{S,t}(x)$ form an orthonormal basis in the space
 $L_2(\mathfrak X^{S,t})$ (this $L_2$ is with respect to the uniform
 measure on $\mathfrak X^{S,t})$.

Let
$$ c^{S,t}_{t+}(i)=\sqrt{\left(1-\frac{i}{t+N}\right)\left(1-\frac{i}{T+N-t-1}\right)}.$$
 Note that $c^{S,t}_{t+}(i)$ does not actually depend on $S$, but it is convenient to use it for the notation.

Finally, denote by $v^{S,t}_{t+}$ the $|\mathfrak X^{S,t}|\times
|\mathfrak X^{t+1,S}|$ matrix given by
$$
 v^{S,t}_{t+}(x,y)=\sum\limits_{i\ge 0} c^{S,t}_{t+}(i) \Psi_i^{S,t}(x)\Psi_i^{S,t+1}(y),\quad x\in \mathfrak
X^{S,t},\quad y\in \mathfrak X^{S,t+1}.
$$

The following proposition was proved in \cite[proof of Proposition
5]{Gor}

\begin{proposition}
\label{Proposition_spectral_t+} Let $X=(x_1<\dots<x_N)$ and
$Y=(y_1<\dots<y_N)$ be elements of $\mathcal X^{S,t}$ and $\mathcal
X^{S,t+1}$, respectively. Then
$$
P^{S,t}_{t+}(X,Y)= \frac{\sqrt{\rho_{S,t+1}(Y)}}
{\sqrt{\rho_{S,t}(X)}}
\frac{\det[v^{S,t}_{t+}(x_i,y_j)]_{i,j=1,\dots,N}}
{\prod\limits_{i=0}^{N-1} c^{S,t}_{t+}(i)}
$$
\end{proposition}

Stochastic matrices $P^{S,t}_{t-}$ admit similar spectral
decomposition with transposed matrices.

Denote
$$
c^{S,t}_{t-}:=c^{S,t-1}_{t+}
$$
and
$$
v^{S,t}_{t-}:=(v^{S,t-1}_{t+})^T
$$
(here $(\quad)^T$ means matrix transposition).

\begin{proposition}
\label{Proposition_spectral_t-}
 Let $X=(x_1<\dots<x_N)$ and
$Y=(y_1<\dots<y_N)$ be elements of $\mathcal X^{S,t}$ and $\mathcal
X^{S,t-1}$, respectively. Then
$$
P^{S,t}_{t-}(X,Y)= \frac{\sqrt{\rho_{S,t-1}(Y) }}
{\sqrt{\rho_{S,t}(X)}}
\frac{\det[v^{S,t}_{t-}(x_i,y_j)]_{i,j=1,\dots,N}}
{\prod\limits_{i=0}^{N-1} c^{S,t}_{t-}(i)}
$$
\end{proposition}
\begin{proof}
Recall that while $P^{S,t}_{t+}$ coincides with transition matrix of
the Markov process $X(t)$, matrix elements of $P^{S,t}_{t-}$ are
co-transition probabilities of Proposition \ref{Proposition_co-tr}.
Consequently,
\begin{multline*}
 P^{S,t}_{t-}(X,Y)=\frac{\rho_{S,t-1}(Y)}{\rho_{S,t}(X)} \cdot
 P^{S,t-1}_{t+}(Y,X)
\\=
\frac{\rho_{S,t-1}(Y)}{\rho_{S,t}(X)} \cdot
\frac{\sqrt{\rho_{S,t}(X)}}{ \sqrt{\rho_{S,t-1}(Y)} }
\frac{\det[v^{S,t-1}_{t+}(y_i,x_j)]_{i,j=1,\dots,N}}
{\prod\limits_{i=0}^{N-1} c^{S,t-1}_{t+}(i)}
\\=
 \frac{\sqrt{\rho_{S,t-1}(Y) }}
{\sqrt{\rho_{S,t}(X)}}
\frac{\det[v^{S,t}_{t-}(y_i,x_j)]_{i,j=1,\dots,N}}
{\prod\limits_{i=0}^{N-1} c^{S,t}_{t-}(i)}
\end{multline*}
\end{proof}

Define
$$ c^{S,t}_{S+}(i)=\sqrt{\left(1-\frac{i}{S+N}\right)\left(1-\frac{i}{T+N-S-1}\right)}$$
(these constants do not acually depend on $t$). Set
$$
 v^{S,t}_{S+}(x,y)=\sum\limits_{i\ge 0} c^{S,t}_{S+}(i) \Psi_i^{S,t}(x)\Psi_i^{S+1,t}(y),\quad x\in \mathfrak
X^{S,t},\quad y\in \mathfrak X^{S+1,t}.
$$

Recall that matrices $P^{S,t}_{t+}$ and $P^{S,t}_{S+}$ are connected
by the involution $t\leftrightarrow S$, thus, Propositions
\ref{Proposition_spectral_t+} and \ref{Proposition_spectral_t-}
imply the following statements.

\begin{proposition}
\label{Proposition_spectral_S+} Let $X=(x_1<\dots<x_N)$ and
$Y=(y_1<\dots<y_N)$ be elements of $\mathcal X^{S,t}$ and $\mathcal
X^{S+1,t}$, respectively. Then
$$
P^{S,t}_{S+}(X,Y)= \frac{\sqrt{\rho_{S+1,t}(Y)}}
{\sqrt{\rho_{S,t}(X)}} \frac{\det[
v^{S,t}_{S+}(x_i,y_j)]_{i,j=1,\dots,N}} {\prod\limits_{i=0}^{N-1}
c^{S,t}_{S+}(i)}.
$$
\end{proposition}

\begin{proposition}
\label{Proposition_spectral_S-} Let $X=(x_1<\dots<x_N)$ and
$Y=(y_1<\dots<y_N)$ be elements of $\mathcal X^{S,t}$ and $\mathcal
X^{S-1,t}$, respectively. Then
$$
P^{S,t}_{S-}(X,Y)= \frac{\sqrt{\rho_{S-1,t}(Y) }}
{\sqrt{\rho_{S,t}(X)}} \frac{\det[
v^{S-1,t}_{S-}(y_i,x_j)]_{i,j=1,\dots,N}} {\prod\limits_{i=0}^{N-1}
c^{S,t}_{S-}(i)}.
$$
\end{proposition}

\subsection{Commutativity}

\begin{theorem}
 The families of stochastic matrices $P^{S,t}_{t\pm}$ and
 $P^{S,t}_{S\pm}$ commute, that is
 $$
   P^{S,t}_{t+}\cdot P^{S,t+1}_{S-}=P^{S,t}_{S-}\cdot
   P^{S-1,t}_{t+},
 $$
 $$
   P^{S,t}_{t-}\cdot P^{S,t-1}_{S-}=P^{S,t}_{S-}\cdot
   P^{S-1,t}_{t-},
 $$
$$
   P^{S,t}_{t+}\cdot P^{S,t+1}_{S+}=P^{S,t}_{S+}\cdot
   P^{S+1,t}_{t+},
 $$
 $$
   P^{S,t}_{t-}\cdot P^{S,t-1}_{S+}=P^{S,t}_{S+}\cdot
   P^{S+1,t}_{t-},
 $$
 for any meaningful values of $S$ and $t$.
\end{theorem}

\begin{proof}
Proofs of all four cases are very similar and we consider only the
first one.
\begin{multline*}
 (P^{S,t}_{t+}\cdot P^{S,t+1}_{S-})(X,Y)=\sum\limits_{Z\in\mathcal
 X^{S,t+1}}P^{S,t}_{t+}(X,Z)\cdot P^{S,t+1}_{S-}(Z,Y)\\
 =\frac{\prod\limits_{i>j}(y_i-y_j)
  \sum\limits_{Z\in\mathcal X^{S,t+1}}  \det[U^{S,t}_{t+}(x_i,z_j)]_{i,j=1,\dots,N}
  \det[U^{S,t+1}_{S-}(z_i,y_j)]_{i,j=1,\dots,N} }
  {(T-t)_{N}\cdot (S)_{N}\cdot\prod\limits_{i>j}(x_i-x_j) }
\end{multline*}
Applying Cauchy-Binet identity we obtain
\begin{multline*}
  \sum\limits_{Z\in\mathcal X^{S,t+1}}  \det[U^{S,t}_{t+}(x_i,z_j)]_{i,j=1,\dots,N}
  \det[U^{S,t+1}_{S-}(z_i,y_j)]_{i,j=1,\dots,N}\\ =
  \det[ (U^{S,t}_{t+}\cdot
  U^{S,t+1}_{S-})(x_i,y_j)]_{i,j=1,\dots,N}.
\end{multline*}
Thus,
$$(P^{S,t}_{t+}\cdot P^{S,t+1}_{S-})(X,Y)=\frac{\prod\limits_{i>j}(y_i-y_j)
\det[ (U^{S,t}_{t+}\cdot  U^{S,t+1}_{S-})(x_i,y_j)]_{i,j=1,\dots,N}}
  {(T-t)_{N}\cdot (S)_{N}\cdot\prod\limits_{i>j}(x_i-x_j)}.
$$
Similarly
$$
(P^{S,t}_{S-}\cdot
P^{S-1,t}_{t+})(X,Y)=\frac{\prod\limits_{i>j}(y_i-y_j) \det[
(U^{S,t}_{S-}\cdot  U^{S-1,t}_{t+})(x_i,y_j)]_{i,j=1,\dots,N}}
  {(T-t)_{N}\cdot (S)_{N}\cdot\prod\limits_{i>j}(x_i-x_j)}.$$
Our problem reduces to verifying the equality
$$
U^{S,t}_{t+}\cdot  U^{S,t+1}_{S-}=U^{S,t}_{S-}\cdot U^{S-1,t}_{t+}.
$$

By straightforward computation one proves that
$$
 U^{S,t}_{t+S-}=U^{S,t}_{t+}\cdot  U^{S,t+1}_{S-}=U^{S,t}_{S-}\cdot
 U^{S-1,t}_{t+},
$$
where
$$U_{t+S-}^{S,t}(x,y)=\begin{cases}
  (N+S-1-x)(N+S-2-x),& \text{ if } y=x+1,\\
  (N+S-1-x)(T-t-S+2x+1),&\text{ if } y=x,\\
  x(T-t-S+x),& \text{ if } y=x-1,\\
  0,&\text{ otherwise.}
\end{cases}
$$
\end{proof}

{\bf Remark.} Another way to prove the commutativity is to use the
spectral decomposition introduced in Section
\ref{Section_Spectral_decomp} and to observe that coefficients
$c^{S,t}_{t\pm}(i)$ do not depend on $S$, while coefficients
$c^{S,t}_{S\pm}(i)$ do not depend on $t$. One computes
\begin{multline*}
\sqrt{\frac{w^{S,t}(x)/w^{S-1,t+1}(y)
}{(t+N)(T+N-t-1)(S+N-1)(T+N-S)}}\cdot U_{t+S-}^{S,t}(x,y)\\=
\sum\limits_{i\ge 0} c^{S,t}_{t+}(i)c^{S,t+1}_{S-}(i)
\Psi_i^{S,t}(x)\Psi_i^{S-1,t+1}(y)\\
=\sum\limits_{i\ge 0} c^{S,t}_{S-}(i)c^{S-1,t}_{t+}(i)
\Psi_i^{S,t}(x)\Psi_i^{S-1,t+1}(y).
\end{multline*}

\section{Markov step $S\mapsto S \pm 1$}

\label{Section_M_steps}

In this section we aim to define two new stochastic matrices
$$P_{S\mapsto S+1}^S(X,Y),\quad X\in\Omega(N,T,S),\quad
Y\in\Omega(N,T,S+1)$$ and
$$P_{S\mapsto S-1}^S(X,Y),\quad X\in\Omega(N,T,S),\quad
Y\in\Omega(N,T,S-1)$$ that preserve the measures $\mu(N,T,S)$. Both
$P_{S\mapsto S+1}^S$ and $P_{S\mapsto S-1}^S$ depend on parameters
$N$ and $T$ but we again omit these indices. In Introduction we
called these matrices $P^{\pm}_{a\times b\times c}$ with $a=T-S$,
$b=S$, $c=N$.

Suppose we are given a sequence $X=(X(0),X(1),\dots,X(T))\in
\Omega(N,T,S)$ (recall that $X(t)\in \mathcal X^{S,t}$). Below we
construct a random sequence $Y=(Y(0),\dots,Y(T))\in\Omega(N,T,S+1)$
and therefore define the transition probability (or, equivalently,
stochastic matrix) $P_{S\mapsto S+1}^S(X,Y)$.

First note that $Y(0)\in \mathcal X^{S+1,0}$ and $|\mathcal
X^{S+1,0}|=1$. Thus, $Y(0)$ is uniquely defined. We will perform a
\emph{sequential update}. Suppose $Y(0),Y(1),\dots,Y(t)$ have been
already defined. Define conditional distribution of $Y(t+1)$ given
$X$, $Y(0),Y(1),\dots,Y(t)$ by
\begin{multline}
\label{Two_kinds}
 {\rm Prob}\{Y(t+1)=Z\}=\frac{P^{S+1,t}_{t+}(Y(t),Z)\cdot
 P^{S+1,t+1}_{S-}(Z,X(t+1))}{
 (P^{S+1,t}_{t+}P^{S+1,t+1}_{S-})(Y(t),X(t+1))}\\
 =\frac{P^{S,t+1}_{S+}(X(t+1),Z)\cdot P^{S+1,t+1}_{t-}(Z,Y(t))
 }{
 (P^{S,t+1}_{S+}P^{S+1,t+1}_{t-})(X(t+1),Y(t))}.
\end{multline}
(The second equality follows from $
\rho_{S+1,t+1}(X)P^{S+1,t+1}_{t-}(X,Y)=\rho_{S+1,t}(Y)P^{S+1,t}_{t+}(Y,X)$.)

This definitions follows the idea of \cite[Section 2.3]{DF}, see
also \cite{BF}.

Observe that $(P^{S+1,t}_{t+}P^{S+1,t+1}_{S-})(Y(t),X(t+1))>0$.
Indeed
\begin{multline}
\label{Tr_nonzero} (P^{S+1,t}_{t+}P^{S+1,t+1}_{S-})(Y(t),X(t+1))=
(P^{S+1,t}_{S-}P^{S,t}_{t+})(Y(t),X(t+1))\\ \ge
P^{S+1,t}_{S-}(Y(t),X(t))\cdot P^{S,t}_{t+}(X(t),X(t+1))>0,
\end{multline} because $Y(t)$ was chosen on the previous step so that
$P^{S+1,t}_{S-}(Y(t),X(t))>0$.

One could say that we choose $Y(t+1)$ using conditional distribution
of the middle point in the successive application of
$P^{S+1,t}_{t+}$ and $P^{S+1,t+1}_{S-}$ (or $P^{S,t+1}_{S+}$ and
$P^{S+1,t+1}_{t-}$ ), provided that we start at $Y(t)$ and finish at
$X(t+1)$ (or start at $X(t+1)$ and finish at $Y(t)$).

After performing $T$ updates we obtain the sequence $Y$.

Equivalently, define $P_{S\mapsto S+1}^S$  by (cf. \cite[Section
2.2]{BF})
$$
 P_{S\mapsto S+1}^S(X,Y)=\begin{cases}
 \prod\limits_{t=0}^{T-1}\dfrac{P^{S+1,t}_{t+}(Y(t),Y(t+1))\cdot
 P^{S+1,t+1}_{S-}(Y(t+1),X(t+1))}{
 (P^{S+1,t}_{t+}P^{S+1,t+1}_{S-})(Y(t),X(t+1))},\\
  \quad \quad \text{if } \prod\limits_{t=0}^{T-1}
  (P^{S+1,t}_{t+}P^{S+1,t+1}_{S-})(Y(t),X(t+1))>0,\\
  0, \text{ otherwise.}
  \end{cases}
$$

\begin{theorem}
\label{Th_Mes_pr} The matrix $P^S_{S\mapsto S+1}$ on
$\Omega(N,T,S)\times\Omega(N,T,S+1)$ is stochastic.
 The transition probabilities $P_{S\mapsto S+1}^S(X,Y)$ preserve the uniform measures
 $\mu(N,T,S)$:
 $$
   \mu(N,T,S+1)(Y)=\sum_{X\in\Omega(N,T,S)}P_{S\mapsto
   S+1}^S(X,Y)\mu(N,T,S)(X).
 $$
\end{theorem}
\begin{proof}
First, let us prove that the matrix $P^S_{S\mapsto S+1}$  is
stochastic, equivalently:
\begin{multline}
\label{StF} 1=\sum\limits_{Y\in\Omega(N,T,S+1)} P^S_{S\mapsto
S+1}(X,Y)
\\=\sum_{Y(0),\dots,Y(T)}\prod\limits_{t=0}^{T-1}\dfrac{P^{S+1,t}_{t+}(Y(t),Y(t+1))\cdot
 P^{S+1,t+1}_{S-}(Y(t+1),X(t+1))}{
 (P^{S+1,t}_{t+}P^{S+1,t+1}_{S-})(Y(t),X(t+1))},
\end{multline}
where the summation goes over all $(Y(0),Y(1),\dots,Y(T))\in
\Omega(N,T,S+1)$ such that
\begin{equation}
\label{St_Ineq}
 \prod\limits_{t=0}^{T-1}
  (P^{S+1,t}_{t+}P^{S+1,t+1}_{S-})(Y(t),X(t+1))>0.
\end{equation}

We write \eqref{StF} in the form
\begin{multline*}
\sum_{Y(0),\dots,Y(T-1)}\prod\limits_{t=0}^{T-2}\dfrac{P^{S+1,t}_{t+}(Y(t),Y(t+1))\cdot
 P^{S+1,t+1}_{S-}(Y(t+1),X(t+1))}{
 (P^{S+1,t}_{t+}P^{S+1,t+1}_{S-})(Y(t),X(t+1))}\\ \times\sum\limits_{Y(T)}\dfrac{P^{S+1,T-1}_{t+}(Y(T-1),Y(T))\cdot
 P^{S+1,T}_{S-}(Y(T),X(T))}{
 (P^{S+1,T-1}_{t+}P^{S+1,T}_{S-})(Y(T-1),X(T))}.
\end{multline*}
Summing over $Y(T)$ we obtain
\begin{multline*}
\sum_{Y(0),\dots,Y(T-2)}\prod\limits_{t=0}^{T-3}\dfrac{P^{S+1,t}_{t+}(Y(t),Y(t+1))\cdot
 P^{S+1,t+1}_{S-}(Y(t+1),X(t+1))}{
 (P^{S+1,t}_{t+}P^{S+1,t+1}_{S-})(Y(t),X(t+1))}\\ \times\sum\limits_{Y(T-1)}\dfrac{P^{S+1,T-2}_{t+}(Y(T-2),Y(T-1))\cdot
 P^{S+1,T-1}_{S-}(Y(T-1),X(T-1))}{
 (P^{S+1,T-2}_{t+}P^{S+1,T-1}_{S-})(Y(T-2),X(T-1))}.
\end{multline*}
Next, we want to sum over $Y(T-1)$. Inequality \eqref{St_Ineq}
implies that the summation goes over $Y(T-1)$ such that
$(P^{S+1,T-1}_{t+}P^{S+1,T}_{S-})(Y(T-1),X(T))>0$. Note that
\begin{multline*}
(P^{S+1,T-1}_{t+}P^{S+1,T}_{S-})(Y(T-1),X(T))=(P^{S+1,T-1}_{S-}P^{S,T-1}_{t+})(Y(T-1),X(T))\\
\ge P^{S+1,T-1}_{S-}(Y(T-1),X(T-1))P^{S,T-1}_{t+}(X(t-1),X(T)).
\end{multline*}

Consequently, if ${(P^{S+1,T-1}_{t+}P^{S+1,T}_{S-})(Y(T-1),X(T))}$
vanishes, then $P^{S+1,T-1}_{S-}(Y(T-1),X(T-1))$ vanishes too. Thus,
we may drop the inequality that restricts the summation and sum over
all possible $Y(T-1)$. After that we sum over $Y(T-2)$, and so on.
After summing over all $Y(t)$, $t=T,T-1,\dots,1$ and noticing that
$Y(0)$ has just one possible value we arrive at \eqref{StF}.

Now we prove that the transition probabilities $P_{S\mapsto
S+1}^S(X,Y)$ preserve the uniform measures
 $\mu(N,T,S)$. It is equivalent to
\begin{multline}
\label{Pr_eq}
\mu(N,T,S+1)(Y)=\sum_{X=(X(0),X(1),\dots,X(T))}P_{S\mapsto
   S+1}^S(X,Y)\mu(N,T,S)(X).
\end{multline}
Since $X(t)$ can be viewed as a Markov chain with time $t$,
$$
 \mu(N,T,S)(X)=m_S^0(X(0))\cdot
 P^{S,0}_{t+}(X(0),X(1))\dots P^{S,T-1}_{t+}(X(T-1),X(T)),
$$
where $m_S^0(X(0))$ is the unique probability measure on singleton
$\mathcal X^{S,0}$.

Thus, the right-hand side of \eqref{Pr_eq} is equal to
\begin{multline}
\label{+1_Measure_Conserving}
 \sum_{X(0),\dots,X(T)}
m_S^0(X(0))\prod\limits_{t=0}^{T-1}P^{S,t}_{t+}(X(t),X(t+1))\\
\times
\prod\limits_{t=0}^{T-1}\dfrac{P^{S+1,t}_{t+}(Y(t),Y(t+1))\cdot
 P^{S+1,t+1}_{S-}(Y(t+1),X(t+1))}{
 (P^{S+1,t}_{t+}P^{S+1,t+1}_{S-})(Y(t),X(t+1))}
\end{multline}
Note that
$$m_S^0(X(0))=m_{S+1}^0(Y(0))=P^{S+1,0}_{S-1}(Y(0),X(0))=P^{S+1,T}_{S-}(Y(T),X(T))=1$$
and write \eqref{+1_Measure_Conserving} in the form
\begin{multline*}
m_{S+1}^0(Y(0))\prod\limits_{t=0}^{T-1}P^{S+1,t}_{t+}(Y(t),Y(t+1))\\
\times
\sum\limits_{X(0),\dots,X(T)}\prod\limits_{t=0}^{T-1}\dfrac{P^{S+1,t}_{S-}(Y(t),X(t))P^{S,t}_{t+}(X(t),X(t+1))}
                              {(P^{S+1,t}_{S-}P^{S,t}_{t+})(Y(t),X(t+1))}.
\end{multline*}
(We used the equality
$P^{S+1,t}_{S-}P^{S,t}_{t+}=P^{S+1,t}_{t+}P^{S+1,t+1}_{S-}$.)
Summing first over $X(0)$, then over $X(1)$, and so on, we get
$$
m_{S+1}^0(Y(0))\prod\limits_{t=0}^{T-1}P^{S+1,t}_{t+}(Y(t),Y(t+1)),$$
which is exactly the distribution $\mu(N,T,S+1)(Y)$.

\end{proof}

Similarly to $P_{S\mapsto S+1}$, one defines a transition matrix
$$
P_{S\mapsto S-1}^S(X,Y),\quad X\in\Omega(N,T,S),\quad
Y\in\Omega(N,T,S-1),
$$
by
$$
 P_{S\mapsto S-1}^S(X,Y)=\begin{cases}
 \prod\limits_{t=0}^{T-1}\dfrac{P^{S-1,t}_{t+}(Y(t),Y(t+1))\cdot
 P^{S-1,t+1}_{S+}(Y(t+1),X(t+1))}{
 (P^{S-1,t}_{t+}P^{S-1,t+1}_{S+})(Y(t),X(t+1))},\\
  \quad \quad \text{if } \prod\limits_{t=0}^{T-1}
  (P^{S-1,t}_{t+}P^{S-1,t+1}_{S+})(Y(t),X(t+1))>0,\\
  0, \text{ otherwise.}
  \end{cases}
$$
Similarly to \eqref{Two_kinds} there is another way to write
$P^S_{S\mapsto S-1}$ because of the equality
\begin{multline*}
 \dfrac{P^{S-1,t}_{t+}(Y(t),Y(t+1))\cdot
 P^{S-1,t+1}_{S+}(Y(t+1),X(t+1))}{
 (P^{S-1,t}_{t+}P^{S-1,t+1}_{S+})(Y(t),X(t+1))}\\=
 \dfrac{
 P^{S,t+1}_{S-}(X(t+1),Y(t+1))\cdot P^{S-1,t+1}_{t-}(Y(t+1),Y(t))}{
 (P^{S,t+1}_{S-}P^{S-1,t+1}_{t-})(X(t+1),Y(t))}
\end{multline*}

Similarly to Theorem \ref{Th_Mes_pr} one proves the following claim.

\begin{theorem}
\label{Th_Mes_pr2}
 The matrix $P^S_{S\mapsto S-1}$ on
$\Omega(N,T,S)\times\Omega(N,T,S-1)$ is stochastic.
 The transition probabilities $P_{S\mapsto S-1}^S(X,Y)$ preserve the uniform measures
 $\mu(N,T,S)$:
 $$
   \mu(N,T,S-1)(Y)=\sum_{X\in\Omega(N,T,S)}P_{S\mapsto
   S-1}^S(X,Y)\mu(N,T,S)(X).
 $$
\end{theorem}

{\bf Remark.} The above construction performs sequential update from
$t=0$ to $t=T$. One can equally well update from $t=T$ to $t=0$ by
suitably modifying the definitions. The resulting Markov chains also
preserve the uniform measure $\mu(N,T,S)$, and they are different
from the Markov chains defined above.

\section{Algorithmic description}
\label{Section_Algorithm} In this section we suggest an algorithmic
description of the Markov chain from the previous section.

Denote by $D(a,b,n)$ the probability distribution on
$\{0,1,\dots,n\}$ given by
\begin{multline}
\label{Distribution}
 {\rm
Prob}(\{k\})=D(a,b,n)\{k\}=\frac{\frac{(a)_k}{(b)_k}}
           {\sum\limits_{j=0,\dots,n}\frac{(a)_j}{(b)_j}
           }\\=\frac{a(a+1)\dots(a+k-1)(b+k)\dots(b+n-1)}{\sum_{j=0}^{n}a(a+1)\dots(a+j-1)(b+j)\dots(b+n-1)}.
\end{multline}

\subsection{ Algorithm for $S\mapsto S+1$ step.}

 Suppose we are
given $X=(X(0),X(1),\dots,X(T))\in\Omega(N,T,S)$. We want to
construct $Y=(Y(0),Y(1),\dots,Y(T))\in\Omega(N,T,S+1)$.

In the first place we note that $Y(0)$ is uniquely defined,
$$
 Y(0)=(0,1,\dots,N-1).
$$
Then we perform $T$ sequential updates, i.e. for $t=0,1,\dots T-1$
we construct $Y(t+1)$ using $Y(t)$ and $X(t+1)$. Let us describe
each step.

Let $Y(t)=(y_1<y_2<\dots<y_N)$ and $X(t+1)=(x_1<x_2<\dots<x_N)$. We
are going to construct $Y(t+1)=(z_1<z_2<\dots<z_N)$.

Recall that
$$z_i\in\mathfrak X^{S+1,t+1}=\{x\in\mathbb Z\mid \max(0,t+S-T+2)\le x \le \min(t+N,S+N)\}.$$

$Y(t)$ and $X(t+1)$ satisfy \eqref{Tr_nonzero}. This implies that
$x_i-y_i$ is equal to either $-1$, $0$ or $1$ for every $i$.

$\bullet$ First, consider all indices $i$ such that $x_i-y_i=1$. For
every such $i$ we set $z_i=x_i$.

$\bullet$ Second, consider all indices $i$ such that $x_i-y_i=-1$
and set $z_i=y_i$.

$\bullet$ Finally, consider all remaining indices, i.e. all $i$ such
that $x_i=y_i$. Divide the corresponding $x_i$'s into blocks of
neighboring integers of distance at least one from each other. Call
such a block a $(k,l)$--block, where $k$ is the smallest number in
the block and $l$ is its size. Thus, we have
$$x_i=y_i=k,\quad x_{i+1}=y_{i+1}=k+1,\quad \dots,\quad x_{i+l-1}=y_{i+l-1}=k+l-1$$
and $$ y_{i-1}<k-1,\quad y_{i+l}>k+l.$$

For each $(k,l)$--block we perform the following procedure: consider
random variable $\xi$ distributed according to
${D(k+T-t-S-1,k+1,l)}$ ($\xi$'s corresponding to different
$(k,l)$--blocks are independent). Set $z_i=x_i$ for the first $\xi$
integers of the block (their coordinates are $k,k+1,\dots, k+\xi-1$)
and set $z_i=x_i+1$ for the rest of the block.

\medskip

At Figure 4 we provide an example of constructing $Y(t+1)$ using
$X(t+1)$ and $Y(t)$: there is only one $(k,l)$--block and it splits
into two groups, here $\xi=2$.

\begin{center}
 {\scalebox{0.5}{\includegraphics{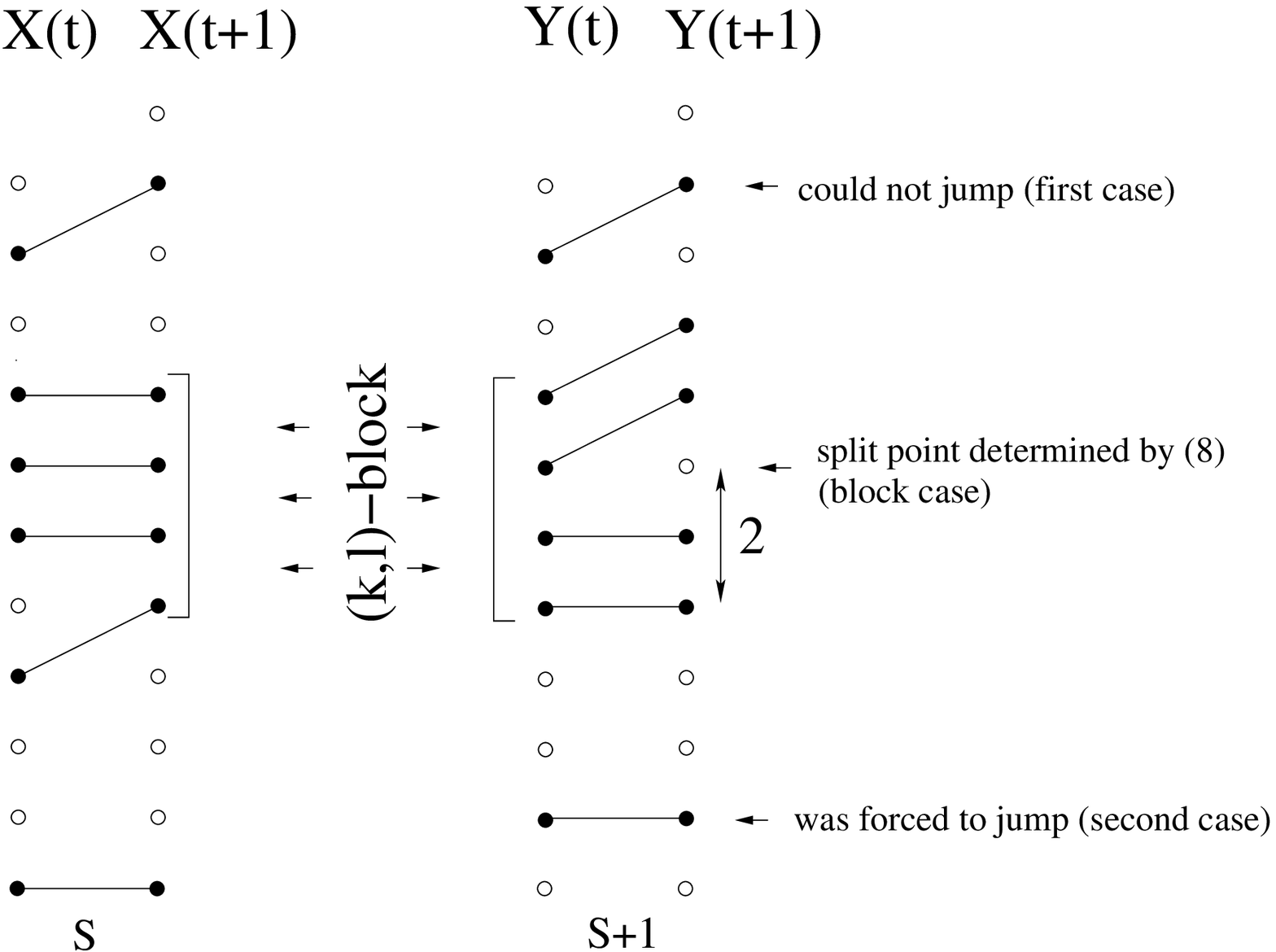}}}

 Figure 4. Example of $(k,l)$--block split, $l=4$, $\xi=2$.
\end{center}

\begin{theorem}
 \label{Alg1}
 The algorithm described above is precisely $S\mapsto S+1$ Markov
 step given by $P^S_{S\mapsto S+1}$.
\end{theorem}

\begin{proof}

As was shown in the previous section, the transition $S\to S+1$
consists of $T$ updates. Namely, given $Y(t)$ and $X(t+1)$ we define
$Y(t+1)$ by
$$
{\rm Prob}\{Y(t+1)=Z\}=\frac{P^{S+1,t}_{t+}(Y(t),Z)\cdot
 P^{S+1,t+1}_{S-}(Z,X(t+1))}{
 (P^{S+1,t}_{t+}P^{S+1,t+1}_{S-})(Y(t),X(t+1))}.
$$
Let $Y(t)=(y_1<y_2<\dots<y_N)$, $X(t+1)=(x_1<x_2<\dots<x_N)$ and we
are defining $Y(t+1)=(z_1<z_2<\dots<z_N)$. Inequality
\eqref{Tr_nonzero} implies that for every $i$ the difference
$x_i-y_i$ is equal to either $-1$, $0$ or $1$. Thus, we have 3
cases:
\begin{enumerate}
\item ${x_i-y_i=1}$: Definitions of $P^{S+1,t}_{t+}$ and $P^{S+1,t+1}_{S-}$ imply that both $z_i-y_i$ and $z_i-x_i$ must be equal to either $0$ or $1$. Consequently,
$z_i$ is uniquely defined and $z_i=x_i$.
\item ${x_i-y_i=-1}$: Again $z_i$ is uniquely defined, and
$z_i=y_i$.
\item ${x_i-y_i=0}$: This is the only nontrivial case. Here one has
two pos\-si\-bi\-li\-ties, either $z_i=x_i$ or $z_i=x_i+1$.

When we pass from $x_i$ to $z_i$, every $(k,l)$--block is split into
at most two groups. Namely, we have $l+1$ possibilities for the
split point $j\in\{0,1,\dots,l\}$ $0,1\dots,l$: $z_i=x_i$ for the
lowest $j$ points of the block and $z_i=x_i+1$ for the other points
of this block.

Now we want to compute the probabilities of different splits of the
blocks. We have
\begin{multline*}
{\rm Prob}\{Y(t+1)=Z\}=\frac{P^{S+1,t}_{t+}(Y(t),Z)\cdot
 P^{S+1,t+1}_{S-}(Z,X(t+1))}{
 (P^{S+1,t}_{t+}P^{S+1,t+1}_{S-})(Y(t),X(t+1))}\\
 =\prod\limits_{i:z_i=y_i}(y_i+T-t-S-1) \cdot
 \prod\limits_{i:z_i=y_i+1}(N+S-y_i)\\ \times
 \prod\limits_{i:z_i=x_i}(N+S-x_i)\cdot \prod\limits_{i:z_i=x_i+1}
 (x_i+1)
 \cdot(\text{factors independent of } Z)
\end{multline*}

This formula implies that blocks split independently. For each
$(k,l)$--block the probability of split position $j$ is equal to
\begin{multline*}
 \prod\limits_{a=k}^{k+j-1}(a+T-t-S-1)(N+S-a)\cdot\prod\limits_{a=k+j}^{k+l-1}(a+1)(N+S-a)
 \\ \times(\text{factors independent of }j)
\end{multline*}
Since $(N+S-a)$ is present in both products, this probability can be
written as
$$
\prod\limits_{a=k}^{k+j-1}(a+T-t-S-1)\cdot\prod\limits_{a=k+j}^{k+l-1}(a+1)\cdot
(\text{factors independent of }j)
$$
which is exactly the distribution $D(k+T-t-S-1,k+1,l)$.
\end{enumerate}
\end{proof}

\subsection{ Algorithm for $S\mapsto S-1$ step}

The $S\mapsto S-1$ step algorithm is very similar to the $S\mapsto
S+1$ one.

Suppose we are given $X=(X(0),X(1),\dots,X(T))\in\Omega(N,T,S)$. We
want to construct $Y=(Y(0),Y(1),\dots,Y(T))\in\Omega(N,T,S-1)$.

As above, note that $Y(0)$ is uniquely defined,
$$
 Y(0)=(0,1,\dots,N-1).
$$
Then we  again perform $T$ sequential updates, i.e. for $t=0,1,\dots
T-1$ we construct $Y(t+1)$ using $Y(t)$ and $X(t+1)$. Let us
describe each step.

Let $Y(t)=(y_1<y_2<\dots<y_N)$ and $X(t+1)=(x_1<x_2<\dots<x_N)$. We
are going to construct $Y(t+1)=(z_1<z_2<\dots<z_N)$.

Recall that
$$z_i\in\mathfrak X^{S-1,t+1}=\{x\in\mathbb Z\mid \max(0,t+S-T)\le x \le \min(t+N,S+N-2)\}.$$

$Y(t)$ and $X(t+1)$ satisfy
$(P^{S-1,t}_{t+}P^{S-1,t+1}_{S+})(Y(t),X(t+1))>0$. This implies that
$x_i-y_i$ is equal to either $0$, $1$ or $2$ for every $i$.

$\bullet$ First, consider all indices $i$ such that $x_i-y_i=0$. For
every such $i$ we set $z_i=x_i$.

$\bullet$ Second, consider all indices $i$ such that $x_i-y_i=2$ and
set $z_i=y_i+1$.

$\bullet$ Finally, consider all remaining indices, i.e. all $i$ such
that $x_i=y_i+1$. Divide the corresponding $x_i$'s into blocks of
neighboring integers of distance at least one from each other. Call
such a block a $(k,l)'$--block, where $k$ is the smallest number in
the block and $l$ is its size. Thus, we have
$$x_i=y_i+1=k,\quad x_{i+1}=y_{i+1}+1=k+1,\quad \dots,\quad x_{i+l-1}=y_{i+l-1}=k+l-1.$$

For each $(k,l)'$--block we perform the following procedure:
consider random variable $\xi$ distributed according to
${D(N+t-k+1,N+S-k-1,l)}$ ($\xi$'s corresponding to different
$(k,l)'$--blocks are independent). Set $z_i=y_i$ for the first $\xi$
integers of the block (their coordinates are $k-1,k,\dots, k+\xi-2$)
and set $z_i=y_i+1$ for the rest of the block.

\begin{theorem}
 The algorithm described above is precisely $S\mapsto S-1$ Markov
 step defined by $P^S_{S\mapsto S-1}$.
\end{theorem}
The proof is similar to Theorem \ref{Alg1} and we omit it.

\subsection{Numeric experiments }
\label{Section_Numeric_exper}

The $S\mapsto S\pm 1$ steps can be used in different ways.

{\bf 1.} Suppose that our aim is to sample a random tiling
(equivalently, random family of paths) $\mathcal T\in\Omega(N,T, S)$
from the uniform measure $\mu(N,T, S)$.

We start from the unique family of paths $\mathcal
T_0\in\Omega(N,T,0)$. Indeed, $|\Omega(N,T,0)|=1$.

Next, we perform $S$ steps. During $r$th step we construct $\mathcal
T_{r}\in\Omega(N,T,r)$ distributed as $\mu(N,T,r)$, using already
constructed family of paths $\mathcal T_{r-1}\in\Omega(N,T,r-1)$.
Theorem \ref{Th_Mes_pr} implies that $\mathcal T_{ S}$ is the
desired random element of $\Omega(N,T, S)$.

Let us estimate the number of operations. Every update takes $O(N)$
operations. For every $S\to S+1$ step we have to perform $T$
updates. Consequently to sample from $\Omega(N,T,S)$ we need
$O(NTS)$ operations.

On Figure 5 we show a random surface generated by our algorithm.
Here $N=1000$, $T=2000$, $S=1000$. It took less than 4 minutes on
our laptop (Intel Core2 Duo 2.2GHz, 2Gb Ram) to generate this
tiling. Theoretically predicted ``arctic ellipse'', see
\cite{CohnLarsenPropp}, is clearly visible on our picture.

\begin{center}
 {\scalebox{0.4}{\includegraphics{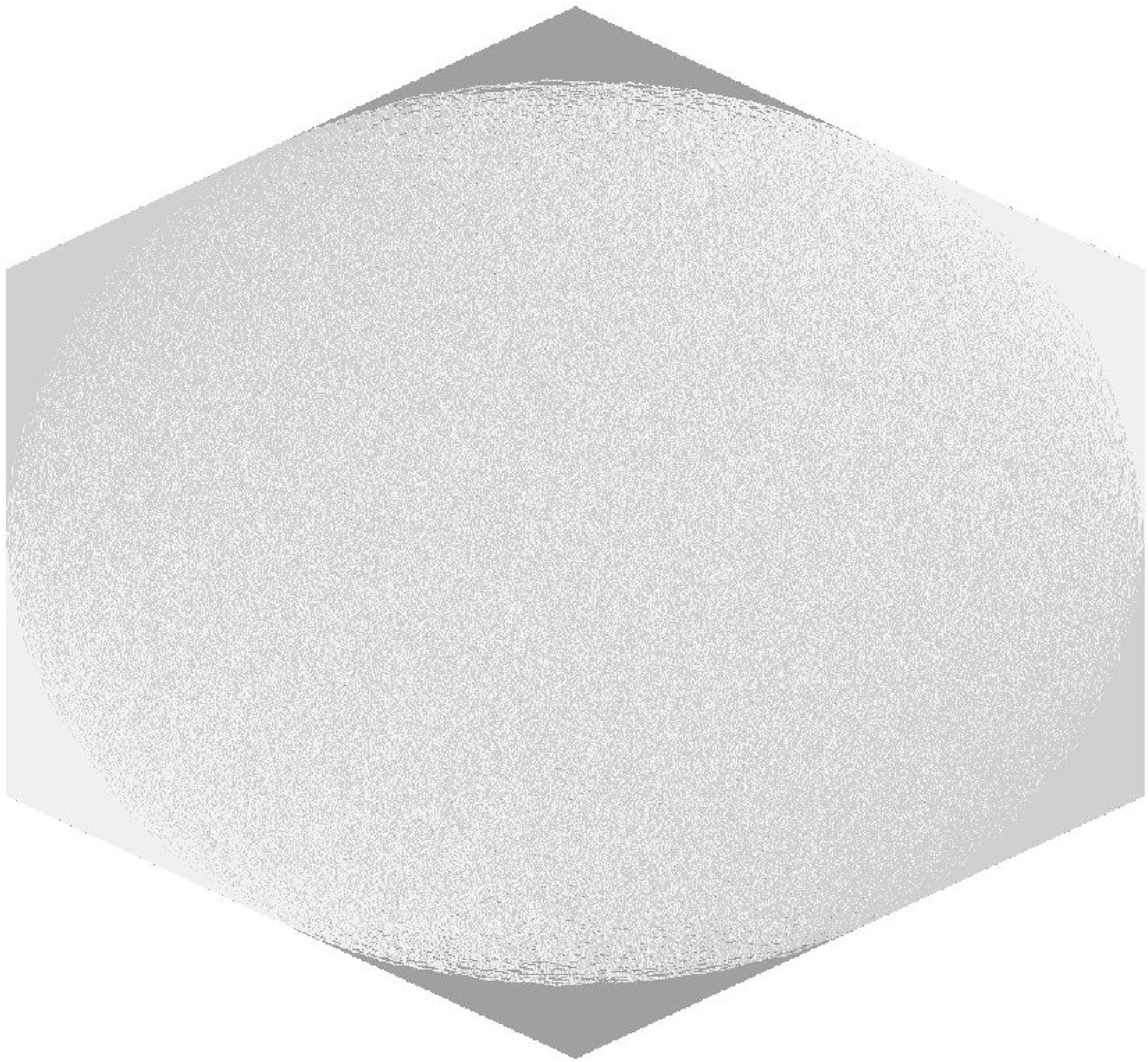}}}

 Figure 5. Random surface corresponding to the big tiling.
\end{center}

{\bf 2.} Using our steps one can also construct equilibrium dynamics
$S\mapsto S+1\mapsto S$ or $S\mapsto S-1\mapsto S$.

On Figure 6 we show the evolution of the ``filled box'' tiling under
$S\mapsto S+1\mapsto S$ dynamics. Here $N=50$, $T=50$, $S=20$.

\begin{tabular}{cc}
 Original tiling& After 20 steps\\
 {\scalebox{0.32}{\includegraphics{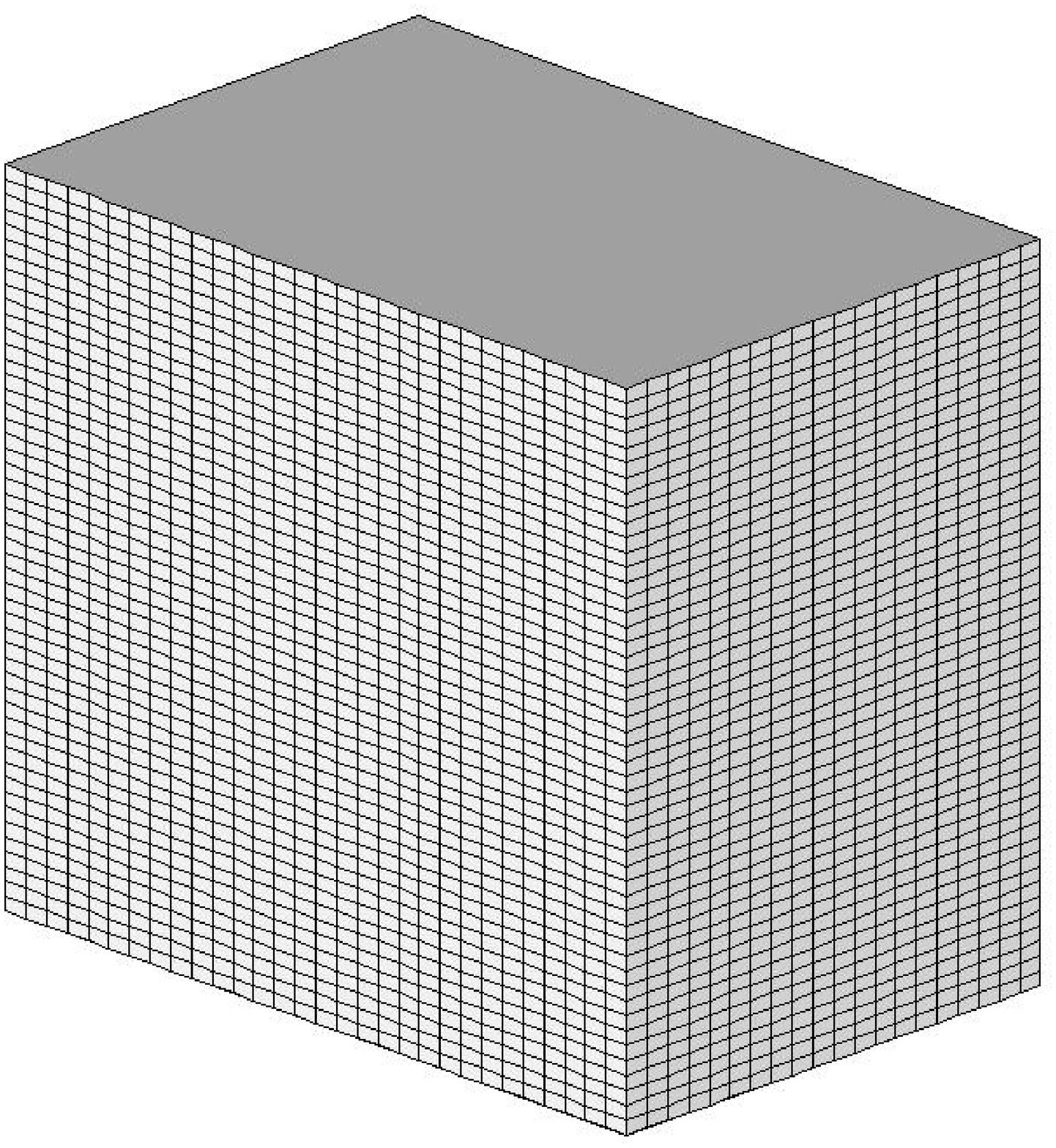}}}&
 {\scalebox{0.32}{\includegraphics{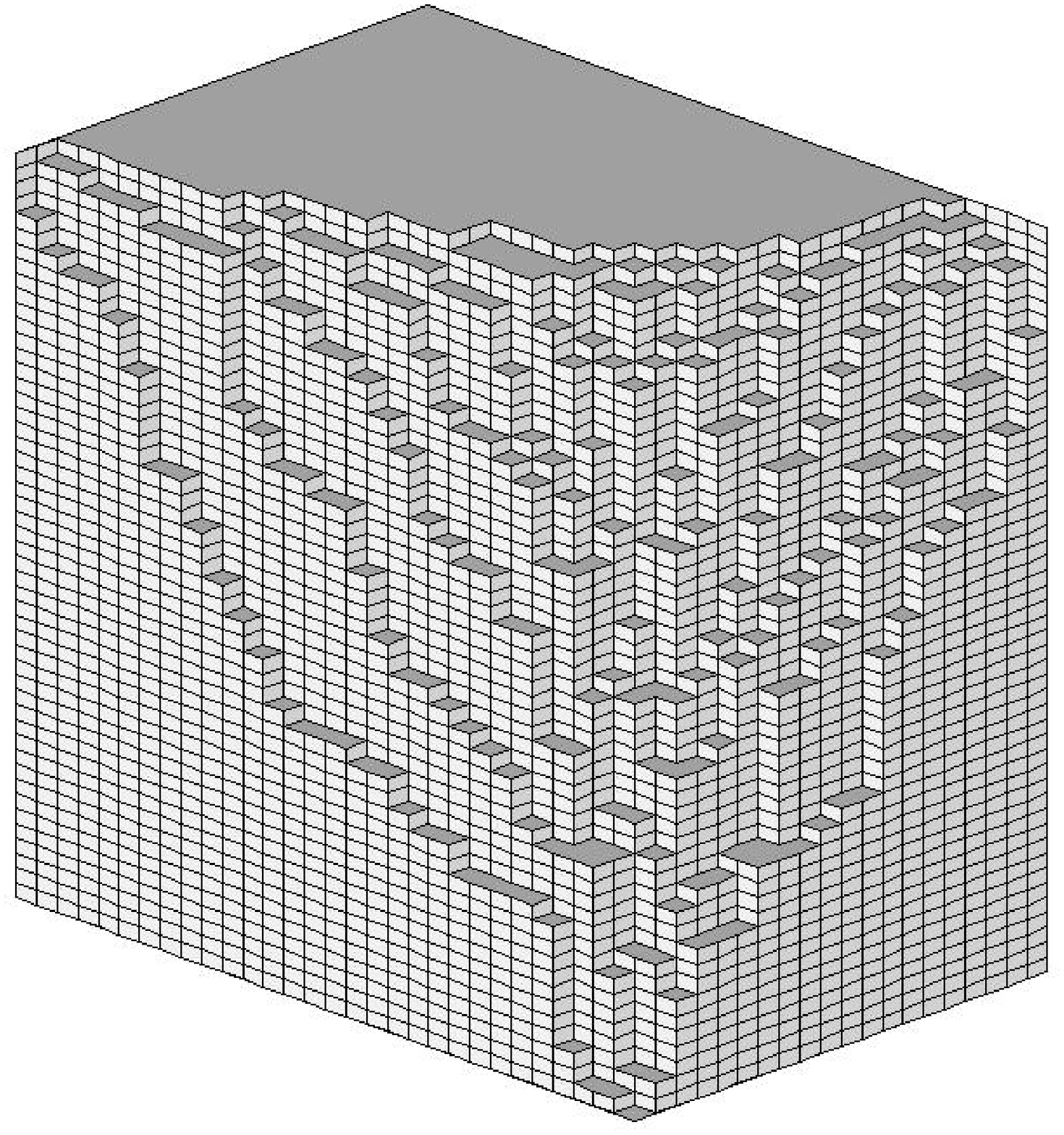}}}\\
 After 100 steps& After 1000 steps\\
  {\scalebox{0.32}{\includegraphics{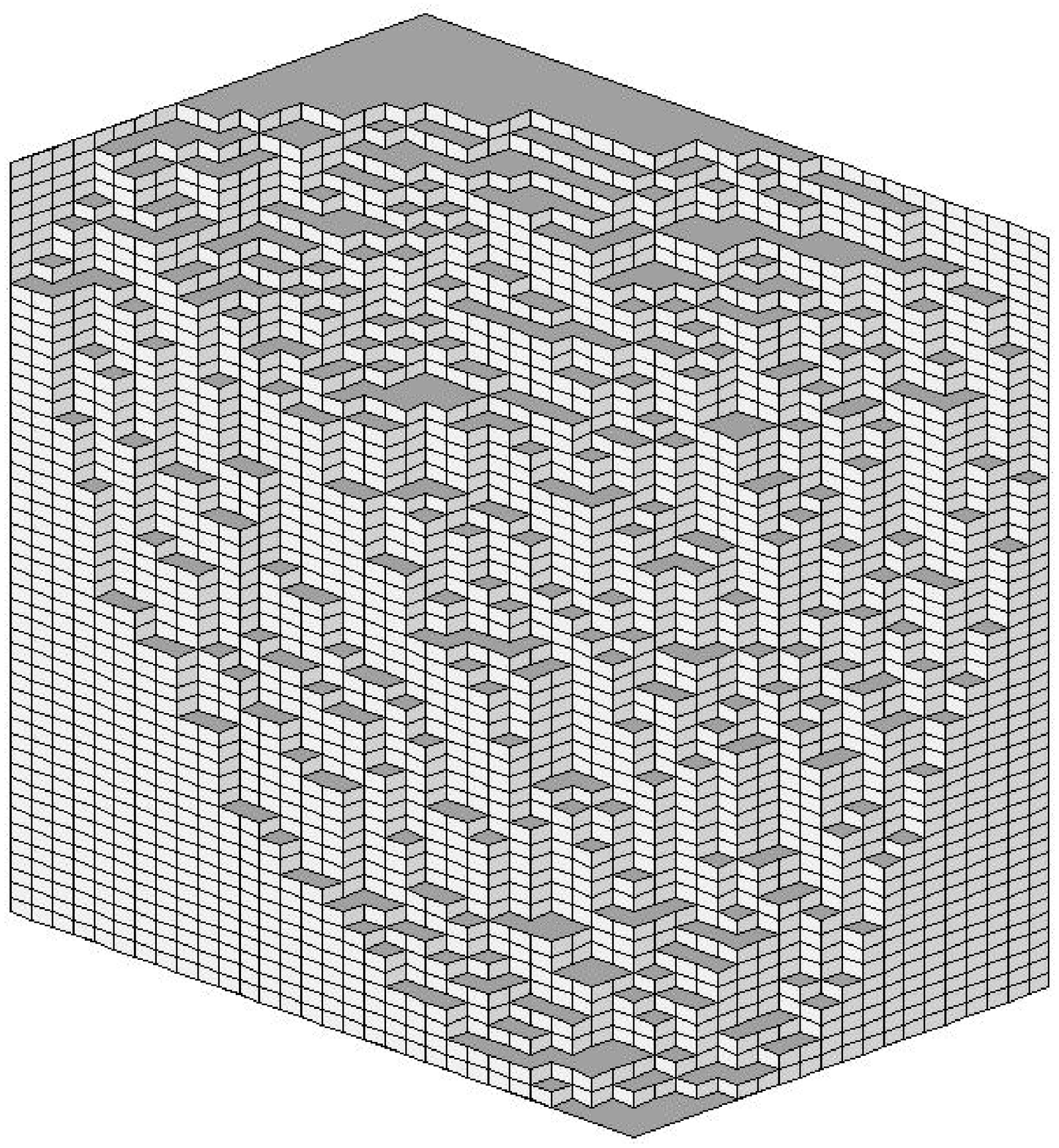}}}&
  {\scalebox{0.32}{\includegraphics{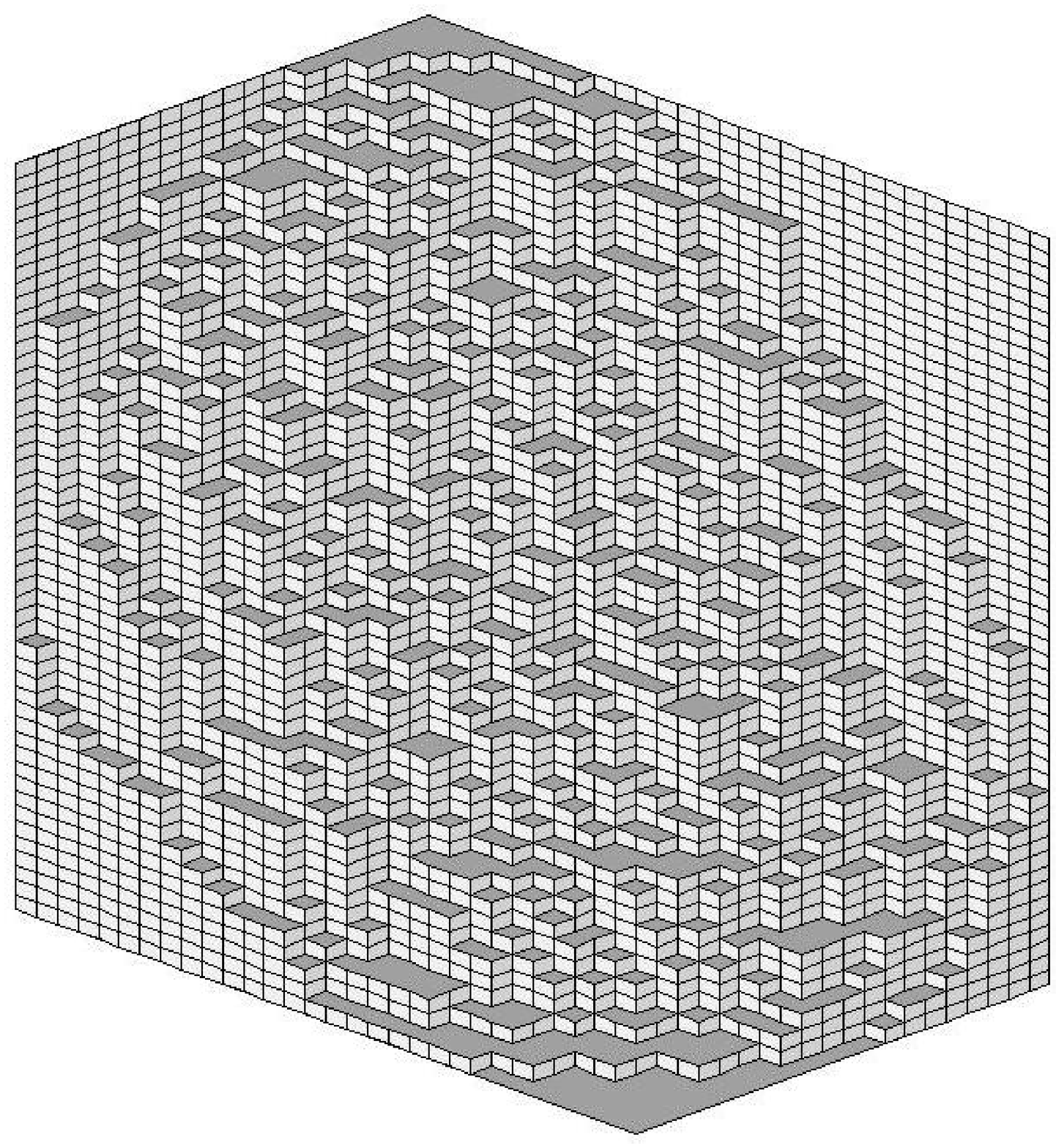}}}
\end{tabular}
\begin{center}
  Figure 6. Evolution of tiling under $S\mapsto S+1\mapsto S$
  dynamics.
\end{center}

\section{General 2-dimensional dynamics, its sections and correlation functions}
\label{Section_General_process}

\subsection{Construction of dynamics and its correlation functions}
In this section we construct a family of Markov chains on the spaces
$\Omega(N,T,S)$ using the transition probabilities $P_{S\mapsto
S+1}^S(X,Y)$ and $P_{S\mapsto S-1}^S(X,Y)$ introduced in Section 3.
Namely, we want to combine $S\mapsto S+1$ and $S\mapsto S-1$ steps.
To fix the order of ``$+1$'' and ``$-1$'' steps we introduce an
auxiliary sequence $\{\epsilon_i\}$ of $+1$'s and $-1$'s.

Formally, let $0\le S_0\le T$ and let
$\epsilon=\{\epsilon_i\}_{i=1,2,\dots}$ be an arbitrary finite or
infinite sequence of $+1$'s and $-1$'s such that
$$0\le S_0+\sum_{n=1}^m \epsilon_i\le T$$ for every $m$. The last
condition is necessary to ensure that the state spaces
$\Omega(N,T,S)$ of the process are not empty.

For any integer $r$ we denote
$$
S(r)=\begin{cases}
       S_0,&r=0,\\
       S_0+\sum_{i=1}^r\epsilon_i,& r>0.
     \end{cases}
$$

Given $S_0$ and $\epsilon_i$  let us define a Markov chain
$$\mathfrak M_{S_0,\epsilon}(r),\quad r=0,1,\dots$$
$\mathfrak M_{S_0,\epsilon}(r)$ takes values in $\Omega(N,T,S(r))$,
its initial distribution is $\mu(N,T,S_0)$:
$$
 {\rm Prob}\{M_{S_0,\epsilon}(0)=X\}=\mu(N,T,S_0)(\{X\}).
$$
Transition probabilities of our process are given by
$$
{\rm Prob}\{\mathfrak M_{S_0,\epsilon}(r+1)=Y\mid\mathfrak
M_{S_0,\epsilon}(r)=X\}=\begin{cases}
            P_{S\mapsto S+1}^{S(r)}(X,Y),&\text{if
            }\epsilon_{r+1}=1\\
            P_{S\mapsto S-1}^{S(r)}(X,Y),&\text{if
            }\epsilon_{r+1}=-1.
          \end{cases}
$$

Theorems \ref{Th_Mes_pr} and \ref{Th_Mes_pr2} imply that one-time
distributions of $\mathfrak M_{S_0,\epsilon}(r)$ are exactly
$\mu(N,T,S(r))$.

{\bf Example 1.} If $S_0=0$ and $\{\epsilon_i\}=\{1,1,1,\dots\}$,
then $\mathfrak M_{S_0,\epsilon}(r)$ is precisely the chain used for
the random tiling sampling in Section \ref{Section_Numeric_exper}.

{\bf Example 2.} If we set $\epsilon_i=(-1)^i$ and then restrict
$\mathfrak M_{S_0,\epsilon}(r)$ on even $r$ we get a stationary
Markov chain from Section \ref{Section_Numeric_exper}.

Recall that $X\in\Omega(N,T,S)$ is a sequence
$X=(X(0),X(1),\dots,X(T))$. Given a trajectory of the Markov chain
$\mathfrak M_{S_0,\epsilon}(r)$ we can construct a point
configuration in $\mathbb Z^3$ with coordinates $(r,t,x)$ such that
the point $(r_0,t_0,x_0)$ is occupied if and only if $x_0\in
\left(\mathfrak M_{S_0,\epsilon}(r_0)\right)(t_0)$. Thus, our Markov
chain defines a measure on such point configuration or,
equivalently, a random point process in $\mathbb Z^3$. Denote it by
$\mathbb M$.

Define the $n$th \emph{correlation function} of $\mathbb M$ by
\begin{multline*}
 R_n(r_1,t_1,x_1;r_2,t_2,x_2;\dots;r_n,t_n,x_n)\\
 ={\mathbb M}\{{\EuScript M}\in{\rm Conf}(\mathbb Z^3)\mid (r_1,t_1,x_1)\in{\EuScript M},
 (r_2,t_2,x_2)\in{\EuScript M},\dots,(r_n,t_n,x_n)\in{\EuScript M}\}.
\end{multline*}

 These correlation functions uniquely define the process $\mathbb
 M$.

 Through the rest of the paper we concentrate on computation of
 correlation functions $R_n$. Unfortunately, we are not able to
 fully describe $R_n$ for all possible arguments. But we
 can compute $R_n$ on certain two-dimensional sections of $\mathbb Z^3$.

 The main result of this section is the following statement.
  \begin{theorem}
  \label{Th_cor_func}
  Let $r_1\le r_2\le \dots \le r_n$, $t_1\ge t_2\ge\dots\ge
  t_n$. Then
$$
   R_n(r_1,t_1,x_1;r_2,t_2,x_2;\dots;r_n,t_n,x_n)=\det\left[K(r_i,t_i,x_i;r_j,t_j,x_j)\right]_{i,j=1,\dots,n},
$$
 where
   \begin{gather*}
K(r,t,x;r',t',x')=\sum_{i=0}^{N-1}\frac{1}{c_i^{r',t';r,t}}\Psi_{i}^{S(r),t}(x)\Psi_{i}^{S(r'),t'}(y)
 ,\\ \text{ if } r\ge r',t\le t';
 \\ K(r,t,x;r',t',x')=-\sum_{i\ge N}c_i^{r,t;r',t'} \Psi_{i}^{S(r),t}(x)\Psi_{i}^{S(r'),t'}(y),
 \\ \text{ if } r<r'\text{ or } r=r', t>t';
 \\ c_i^{r,t;r,t}=1,
 \\
 c_i^{r,t;r',t'}=\prod\limits_{k=t'+1}^{t}c_{t-}^{S(r),k}(i)\prod\limits_{k=r}^{r'-1}c_{S\epsilon_{k+1}}^{S(k),t}(i),
 \end{gather*}
 where $S\epsilon_{k+1}$ stands for $S+$ if $\epsilon_{k+1}=+1$ and $S-$
 otherwise.
  \end{theorem}

In Section \ref{Section_Bulk} we will study the bulk asymptotics of
these correlation functions.

For $r_1=r_2=\dots=r_n$, Theorem \ref{Th_cor_func} was obtained in
\cite{J}, \cite{JN}, \cite{Gor}.

\subsection{Admissible sections}
We call a sequence $\mathcal A=( (r_0,t_0),(r_1,t_1),\dots (r_n,t_n)
)$ an \emph{admissible section} of $\mathbb Z^2$ provided that
\begin{enumerate}
\item $(r_i,t_i)\in\{0,1\dots\}\times\{0,1,\dots,T\}$
\item $r_0\le
r_1\le\dots\le r_n$
\item $t_0\ge
t_1\ge\dots\ge t_n$
\item For every $i=0,1,\dots,n-1$ either $r_{i+1}=r_i+1$, $t_{i+1}=t_i$ or
$r_{i+1}=r_i$, $t_{i+1}=t_i-1$
\end{enumerate}

Figure 6 gives an example of an admissible section.

\begin{center}
 {\scalebox{0.5}{\includegraphics{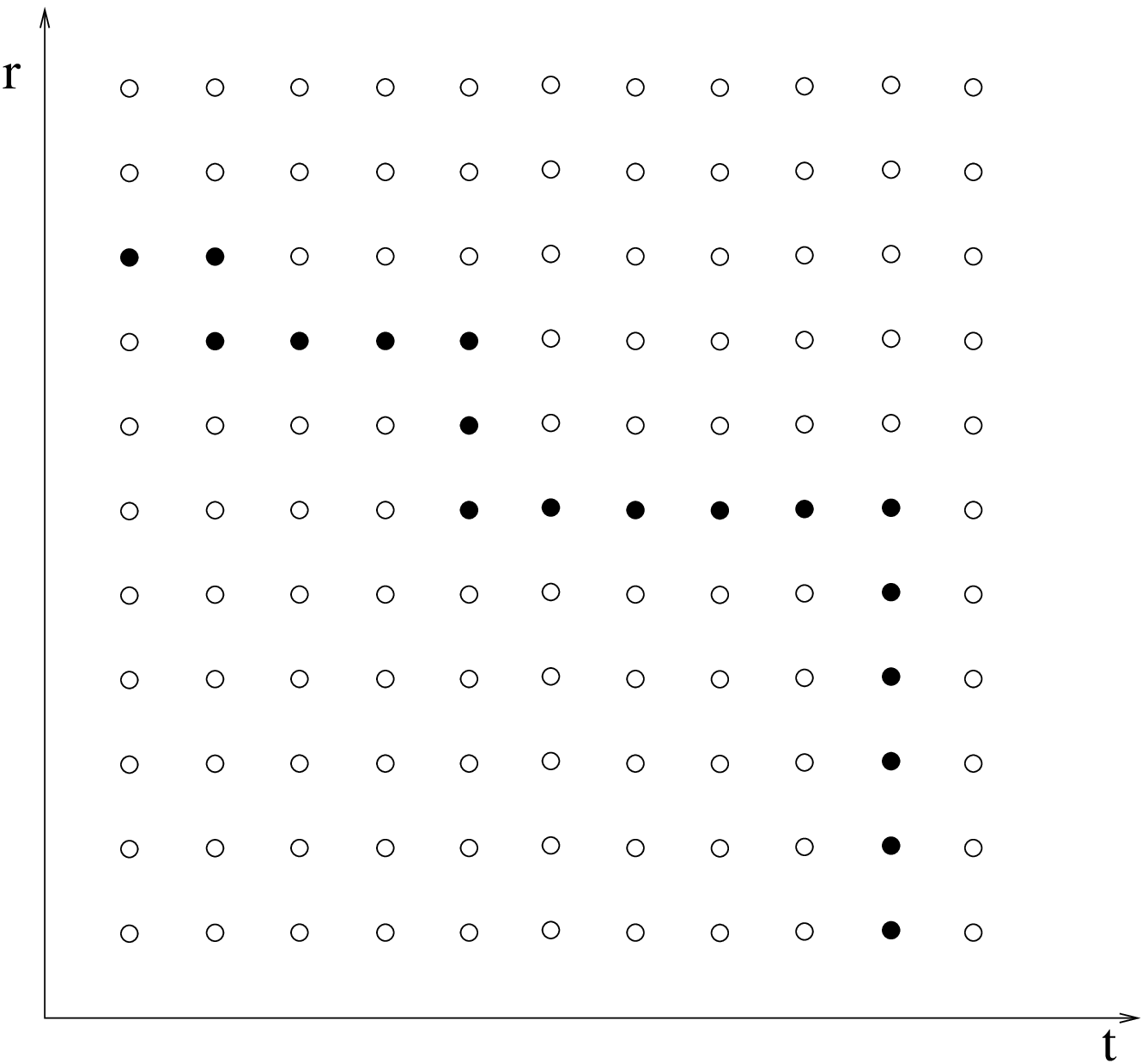}}}

 Figure 6. Admissible section.
\end{center}

Given an admissible section $\mathcal A$ we introduce a
\emph{sectional process} $\EuScript M_{ \mathcal A}(h),\quad
h=0,1\dots,n$, by
$$
 \EuScript M_{\mathcal
  A}(h)=\left(\mathfrak M_{S_0,\epsilon}(r_h)\right)(t_h)
$$

\begin{theorem}
 $\EuScript M_{ \mathcal
  A}(h)$ is a Markov chain
 with initial distribution
 $$
  {\rm Prob}\{\EuScript M_{\mathcal
  A}(0)=X\}=\rho_{N,T}^{S(r_0),t_0}(X)
 $$
 and transition probabilities
 \begin{multline*}
  {\rm Prob}\{\EuScript M_{\mathcal
  A}(h+1)=Y\mid  \EuScript M_{\mathcal
  A}(h)=X\}\\=\begin{cases}
   P^{S(r_h),t_h}_{t-}(X,Y),& \text{ if }t_{h+1}=t_h-1,\\
   P^{S(r_h),t_h}_{S+}(X,Y),& \text{ if }S(r_{h+1})=S(r_h)+1,\\
   P^{S(r_h),t_h}_{S-}(X,Y),& \text{ if }S(r_{h+1})=S(r_h)-1.
  \end{cases}
 \end{multline*}
\end{theorem}
\begin{proof}
This theorem follows from \cite[Proposition 2.7]{BF}. Let us explain
the correspondence between our notations and notations of \cite{BF}.
To avoid confusions we denote by $\tau$ the time variable $t$ from
\cite{BF}. Then $\tau$ corresponds to time $r$ of Markov chain
$\mathfrak M_{S_0,\epsilon}(r)$; $k$ of \cite{BF} corresponds to
$t$; state space ${\cal S}_k(\tau)$ is $\mathcal X^{S(\tau),k}$;
matrices $\Lambda_{k-1}^k(\,\cdot\,,\,\cdot\,\mid \tau)$ of
\cite{BF} are $P^{S(\tau),k}_{t-}(\,\cdot\,,\,\cdot\,)$; matrices
$P_k(\,\cdot\,,\,\cdot\,\mid \tau)$ are
$P^{S(\tau),k}_{S\epsilon_{\tau+1}}(\,\cdot\,,\,\cdot\,)$; the
commutation relations $ \Lambda^k_{k-1}P_{k-1}=P_{k}\Lambda^k_{k-1}
$ are exactly
$P^{S,t}_{S\epsilon}P^{S+\epsilon,t}_{t-}=P^{S,t}_{t-}P^{S,t-1}_{S\epsilon}$
above; $\tau_0$ is $r_0$ and $m_n(\,\cdot\,)$ of \cite{BF}
corresponds to $\rho_{S(r_0),t_0}$.
\end{proof}

\subsection{Correlation functions}
To compute the correlation functions $R_n$ we are going to use a
variant of the Eynard-Metha theorem (see \cite{EM} and \cite[Section
7.4]{BO}). Let us state it first.

\begin{proposition}
\label{proposition_EM} Assume that for every time moment $h$ we are
given
 an orthonormal system $\{g^h_n\}_{n\ge 0}$ in linear space
$l_2(\{0,1, \dots,L\})$  and a set of numbers
 $c_0^h,c_1^h, \dots$. Denote
$$
v_{h,h+1}(x,y)=\sum_{n\ge 0}c_n^h g^h_n(x)g^{h+1}_n(y).
$$

 Assume also that we are given a discrete time  Markov process
${\cal P}_h$ taking values in $N$--tuples of elements of the set
$\{0,1, \dots,L\}$, with one--dimensional distributions
$$\Bigl(\det\left[g_{i-1}^h(x_j)\right]_{i,j=1,\dots ,N}\Bigr)^2$$
 and transition probabilities

 $$\frac{\det\left[v_{h,h+1}(x_i,y_j)\right]_{i,j=1,\dots ,N}
 \det\left[g_{i-1}^{h+1}(y_j)\right]_{i,j=1,\dots ,N}} {\det\left[g_{i-1}^h(x_j)\right]_{i,j=1,\dots ,N}
  \prod\limits_{n=0}^{N-1}c_n^h}. $$

  Then
$$
 {\rm Prob}\{x_1\in {\cal P}_{k_1},\dots, x_n\in {\cal P}_{k_n}\}\\
 =\det\left[K(k_i,x_i;k_j,x_j)\right]_{i,j=1,\dots ,n},
$$ where
\begin{gather*}
K(k,x;l,y)=\sum_{i=0}^{N-1}\frac{1}{c_i^{l,k}}g_i^k(x)g_i^l(y)
 ,\, k\ge l; \\ K(k,x;l,y)=-\sum_{i\ge N}c_i^{k,l} g_i^k(x)g_i^l(y),
 \,k<l; \\ c_i^{k,k}=1,\, c_i^{k,l}=c_i^k\cdot c_i^{k+1}\cdot\dots\cdot
 c_i^{l-1}.
 \end{gather*}
 \end{proposition}

Now set ${\cal P}_h:=\EuScript M_{\mathcal
  A}(h)$. Then we can take orthonormal functions $g^h_n$ to be the
  functions
  $\Psi_n^{S(r_h),t_h}(x)$ defined in Section \ref{Section_Spectral_decomp}, and
\begin{equation}
  \label{c_def}
   c_i^h:=\begin{cases}
   c_{t-}^{S(r_h),t_h}(i),& \text{ if }t_{h+1}=t_h-1,\\
   c_{S+}^{S(r_h),t_h}(i),& \text{ if }S(r_{h+1})=S(r_h)+1,\\
   c_{S-}^{S(r_h),t_h}(i),& \text{ if }S(r_{h+1})=S(r_h)-1,
  \end{cases}
\end{equation}
$$
   v_{h,h+1}(x,y)
   :=\begin{cases}
   v_{t-}^{S(r_h),t_h}(x,y),& \text{ if }t_{h+1}=t_h-1,\\
   v_{S+}^{S(r_h),t_h}(x,y),& \text{ if }S(r_{h+1})=S(r_h)+1,\\
   v_{S-}^{S(r_h),t_h}(y,x),& \text{ if }S(r_{h+1})=S(r_h)-1.
  \end{cases}
$$

  \begin{proposition}
   Markov chain $\EuScript M_{\mathcal A}(h)$ satisfies the assumptions of Proposition \ref{proposition_EM}.
  \end{proposition}
  \begin{proof}
  By Theorem \ref{Th_4stochmatr} one-dimensional distributions of $\EuScript M_{\mathcal
  A}(h)$ are given by the measures $\rho_{S(r_h),t_h}$.
  These are Hahn \emph{orthogonal polynomial ensembles} (see e.g.
  \cite{J},\cite{Gor}). It is well known (see for example
  \cite[Section 2.7]{Konig}) that
  this distribution  can be written in the form
  $$\rho_{S(r_h),t_h}(X)=\left(\det\left[\Psi_{i-1}^{S(r_h),t_h}(x_j)\right]_{i,j=1,\dots ,N}\right)^2.$$
  Propositions \ref{Proposition_spectral_t-}--\ref{Proposition_spectral_S-} imply that the transition probabilities
  can be expressed in the required form too.
  \end{proof}

  Applying Proposition \ref{proposition_EM} we obtain the following
  \begin{proposition}
  \label{Proposition_Cor_func_simple}
   $$
   {\rm Prob}\{x_1\in \EuScript M_{\mathcal A}(h_1),\dots,x_n\in \EuScript M_{\mathcal
   A}(h_n)\}=\det\left[K(k_i,x_i;k_j,x_j)\right]_{i,j=1,\dots,n},
   $$
   \begin{gather*}
K(k,x;l,y)=\sum_{i=0}^{N-1}\frac{1}{c_i^{l,k}}\Psi_{i}^{S(r_k),t_k}(x)\Psi_{i}^{S(r_l),t_l}(y)
 ,\, k\ge l; \\
K(k,x;l,y)=-\sum_{i\ge N}c_i^{k,l}
\Psi_{i}^{S(r_k),t_k}(x)\Psi_{i}^{S(r_l),t_l}(y),
 \,k<l; \\ c_i^{k,k}=1,\, c_i^{k,l}=c_i^k\cdot c_i^{k+1}\cdot\dots\cdot
 c_i^{l-1},
 \end{gather*}
 and coefficients $c_i^{j}$ are given by \eqref{c_def}.
  \end{proposition}

 \begin{proof}[Proof of Theorem \ref{Th_cor_func}]
  If $r_1\le r_2\le \dots \le r_n$, $t_1\ge t_2\ge\dots\ge t_n$, then the sequence
  $\{(r_1,t_1),\dots,(r_n,t_n)\}$ can be
  included into some admissible section $\mathcal A$. Applying
  Proposition \ref{Proposition_Cor_func_simple} and substituting the values of all parameters $c_i^j$ we obtain the desired
  formula for correlation functions $R_n$.
  \end{proof}

\section{Bulk limits}
\label{Section_Bulk}

In this section we aim to compute so-called ``bulk limits'' of the
correlation functions introduced in Section
\ref{Section_General_process}.

Note that while in the previous sections parameters $N,T$ were
always the same, through this section they will change.

We are interested in the following limit regime: Let us fix positive
numbers $\tilde S_0$, $\tilde T$, $\tilde N$, $\tilde t$, $\tilde
x$.
 Introduce a small parameter $\varepsilon \ll 1$, and set
 $$
   S_0= \tilde S_0\varepsilon^{-1} + o(\varepsilon^{-1}),\quad
   T= \tilde T\varepsilon^{-1} + o(\varepsilon^{-1}),\quad
   N= \tilde N\varepsilon^{-1} + o(\varepsilon^{-1}).
 $$
  Consider also integer valued functions $t_i=t_i(\varepsilon)$ and
  $x_i=x_i(\varepsilon)$, $i=1,\dots,n$, such that
  $$
   \lim\limits_{\varepsilon\to 0} \varepsilon
   t_i(\varepsilon)=\tilde t,\quad \lim\limits_{\varepsilon\to 0} \varepsilon x_i(\varepsilon)=\tilde x,\quad
   i=1,\dots,n,
  $$
  and pairwise differences $t_i-t_j$, and $x_i-x_j$ do not depend on
  $\varepsilon$.

  Then correlation functions $R_n$ defined in Theorem \ref{Th_cor_func} tend to a limit
 $\hat R_n$ which
 depends on the parameters of the limit regime $\tilde S_0$, $\tilde T$, $\tilde N$, $\tilde t$, $\tilde
 x$.

 We consider the region where the limit correlation functions are nontrivial. This region is commonly referred to as ``bulk''
 and it is characterized (see e.g. \cite{Gor}), by the fact that  the expression
  $$\frac{-\tilde N(\tilde N+\tilde T)+(-\tilde x+\tilde S_0+\tilde N)
   (\tilde t+\tilde N-\tilde x)+
           \tilde x(\tilde T+\tilde x-\tilde S_0-\tilde t)}
   {2\sqrt{\tilde x(-\tilde x+\tilde S_0+\tilde N)
      (\tilde t+\tilde N-\tilde x)(\tilde x+\tilde T-\tilde S_0-\tilde t)}}.
 $$
 is strictly between $-1$ and $1$. Denote by $\phi=\phi(\tilde S_0,\tilde
 T,\tilde N,\tilde t,\tilde x)$ the arccosine of this expression.

 Let us denote by $Q_i$ the triplet $(r_i,t_i,x_i)$.

\begin{theorem}
\label{Th_bulk_limit} Let  $r_1\le r_2\le \dots \le r_n$, $t_1\ge
t_2\ge\dots\ge
  t_n$ and $R_n(Q_1,\dots,Q_n)$ be defined as in Theorem \ref{Th_cor_func}.
  Then
  $$
  \lim\limits_{\varepsilon\to 0}R_n(Q_1,\dots,Q_n)=\det[K_{ij}^\bu]_{i,j=1,\dots,n},
  $$
  where for $(i,j)$ such that $r_i<r_j$, or $r_i=r_j$, $t_i>t_j$
$$ K^\bu_{ij}= \frac{1}{2\pi i}\int_{\gamma_-}
 \left(1+c_1 z\right)^{t_i-t_j}
\cdot\prod\limits_{k=r_i+1}^{r_j} \left(1+c_2 z^{-\epsilon_k}\right)
  \cdot
  \frac{dz}{z^{x_i-x_j+1}}
$$
  and for $(i,j)$ such that $r_i\ge r_j$, $t_i \le t_j$
   $$K^\bu_{ij}=
\frac{1}{2\pi i}\int_{\gamma_+}
 \left(1+c_1 z\right)^{t_i-t_j}
\cdot\prod\limits_{k=r_j+1}^{r_i} \left(1+c_2
z^{-\epsilon_k}\right)^{-1}
  \cdot \frac{dz}{z^{x_i-x_j+1}},
$$

$$
 c_1=\sqrt{\frac{\tilde x(\tilde S_0+\tilde N-\tilde x)}{(\tilde T-\tilde t-\tilde S_0+\tilde x)(\tilde t+\tilde N-\tilde
x)}},
$$

$$
c_2=\sqrt{\frac{\tilde x(\tilde t+\tilde N-\tilde x)}{(\tilde
T-\tilde t-\tilde S_0+\tilde x)(\tilde S_0+\tilde N-\tilde x)}}.
$$

Here $\gamma_{\pm}$ are contours in $\mathbb C$ joining
$e^{-i\phi(\tilde S_0,\tilde
 T,\tilde N,\tilde t,\tilde x)}$ and $e^{i\phi(\tilde S_0,\tilde
 T,\tilde N,\tilde t,\tilde x)}$ and
crossing $\mathbb R_{\pm}$, respectively.
\end{theorem}

{\bf Comments.} The limiting correlation functions
$$\hat
  R_n(Q_1,\dots,Q_n)=
  \det[K^\bu_{ij}]_{i,j=1,\dots,n}$$
  are correlation
functions of the limit process defined on a fixed admissible
section. The proof of the existence of the limit process can be
found for instance in \cite[Lemma 4.1]{Bor}. The limit process
satisfies certain Gibbs property, see \cite{BS}.

The case $r_1=r_2=\dots=r_n=0$ was thoroughly studied in \cite{Gor}.
The limit process for this case first appeared in \cite{OR}. Our
argument is based on the facts proved in \cite{Gor}.
\begin{proof}
Our first goal is to find the limit of the correlation kernels

$$\lim\limits_{\varepsilon\to 0}
K(r_i,t_i(\varepsilon),x_i(\varepsilon);
r_j,t_j(\varepsilon),x_j(\varepsilon)).
$$

Let us introduce six auxiliary families of $\mathbb Z\times \mathbb
Z$ matrices or, equivalently, operators acting in $l_2(\mathbb Z)$.
Set
$$
{\mathbf P}^{N,T,S,t}(x,y)= \begin{cases} \sum\limits_{k=0}^{N-1}\Psi_k^{S,t}(x)\Psi_k^{S,t}(y), &x,y\in \mathfrak X^{S,t},\\
                 0 &\text{for other } x,y,
   \end{cases}
$$
(recall that $\Psi_k^{S,t}(x)$ and $\mathfrak X^{S,t}$ depend on $N$
and $T$, although these indices are omitted)
$$
 {\mathbf P'}^{N,T,S,t}(x,y)= \begin{cases} -\sum\limits_{k \ge N}\Psi_k^{S,t}(x)\Psi_k^{S,t}(y), &x,y\in \mathfrak X^{S,t},\\
                 0 &\text{for other } x,y.
   \end{cases}
$$

Observe that ${\mathbf P}-{\mathbf P}'={\rm Id}_{\mathfrak
X^{S,t}}$. Define also

$$
   V^{N,T,S,t}_{t+}(x,y)=\begin{cases} \sum\limits_{k \ge 0}c_{t+}^{S,t}(k)\Psi_k^{S,t}(x)\Psi_k^{S,t+1}(y)
        &x\in{\mathfrak X^{S,t}},y\in{\mathfrak X^{S,t+1}},\\
            0 &\text{for other } x,y,
  \end{cases}
 $$
 $$
   V^{N,T,S,t}_{t+}(x,y)=\begin{cases} \sum\limits_{k \ge 0}c_{t-}^{S,t}(k)\Psi_k^{S,t}(x)\Psi_k^{S,t-1}(y)
        &x\in{\mathfrak X^{S,t}},y\in{\mathfrak X^{S,t-1}},\\
            0 &\text{for other } x,y,
  \end{cases}
 $$
 $$
   V^{N,T,S,t}_{S+}(x,y)=\begin{cases} \sum\limits_{k \ge 0}c_{S+}^{S,t}(k)\Psi_k^{S,t}(x)\Psi_k^{S+1,t}(y)
        &x\in{\mathfrak X^{S,t}},y\in{\mathfrak X^{S+1,t}},\\
            0 &\text{for other } x,y,
  \end{cases}
 $$
 $$
   V^{N,T,S,t}_{S-}(x,y)=\begin{cases} \sum\limits_{k \ge 0}c_{S-}^{S,t}(k)\Psi_k^{S,t}(x)\Psi_k^{S-1,t}(y)
        &x\in{\mathfrak X^{S,t}},y\in{\mathfrak X^{S-1,t}},\\
            0 &\text{for other } x,y.
  \end{cases}
 $$

 These are (trivial) extensions of finite matrices $v^{S,t}_{t\pm}$,
 $v^{S,t}_{S\pm}$ of Section \ref{Section_Spectral_decomp}.

 Now let
 $$
 \begin{array}{llll}
   S= \tilde S_0 \varepsilon^{-1} + o(\varepsilon^{-1}),&
   T= \tilde T \varepsilon^{-1}+ o(\varepsilon^{-1}),&
   N= \tilde N \varepsilon^{-1}+ o(\varepsilon^{-1}),\\
   t= \tilde t \varepsilon^{-1}+ o(\varepsilon^{-1}),&
   x= [\tilde x \varepsilon^{-1}]+ \varkappa,&
   y= [\tilde x \varepsilon^{-1}]+ \nu,
  \end{array}
  $$
 and send $\varepsilon$ to $0$. All 6 families tend to some limits. Let us
 denote these limit operators by ${\mathbf{ \hat  P}}(\varkappa, \nu)$, ${\mathbf{ \hat
  P'}}(\varkappa, \nu)$, $\hat V_{t+}(\varkappa, \nu)$, $\hat V_{t-}(\varkappa, \nu)$, $\hat V_{S+}(\varkappa, \nu)$ and $\hat
 V_{S-}(\varkappa, \nu)$, respectively. These limit operators depend on the
 parameters of limit regime $\tilde N$,$\tilde T$,$\tilde S_0$,
 $\tilde t$, $\tilde x$.

As we are dealing with linear operators in $l_2({\mathbb Z})$, it is
convenient  to  employ the Fourier transform
$$F:l_2(\mathbb Z)\mapsto L_2(S^1),$$ where $S^1=\{z\in{\mathbb C}: |z|=1\}$. We denote the images of our operators
under Fourier transform by $F\hat {\mathbf P}$, $F\hat {\mathbf P}'$
and so on.
\begin{proposition}
 When $\varepsilon\to 0$ operators ${\mathbf P}^{N,T,S,t}$ and ${\mathbf
 {P'}}^{N,T,S,t}$ strongly converge to limits ${\mathbf {\hat P}}$ and ${\mathbf {\hat P'}}$,
 respectively.

 $F{\mathbf{ \hat P}}$  is the operator of multiplication by the characteristic function
 of the right arc of the unit circle contained between the angles
 $-\phi$ and $\phi$.

 $F{\mathbf {\hat P'}}$  is the operator of multiplication by the characteristic function
 of the left arc of the unit circle contained between the angles
 $-\phi$ and $\phi$.
\end{proposition}
\begin{proof}
See Section 3.2 in \cite{Gor}.
\end{proof}

\begin{proposition}
\label{Proposition_V_t+_conv} When $\varepsilon\to 0$ operators
$V^{N,T,S,t}_{t+}$ strongly converge to a limit $\hat V_{t+}$.
$F\hat V_{t+}$ is the operator of multiplication by the function
$$\sqrt
     \frac{(\tilde T-\tilde t-\tilde S_0+\tilde x)(\tilde t+\tilde N-\tilde x)}{(\tilde t+\tilde N)(\tilde T+\tilde N-\tilde t)}
            +z^{-1}\cdot\sqrt\frac{(\tilde S_0+\tilde N-\tilde x)\tilde x}{(\tilde t+\tilde N)(\tilde T+\tilde N-\tilde
            t)}.$$
\end{proposition}
\begin{proof}
See Section 3.3 in \cite{Gor}
\end{proof}

\begin{proposition}
When $\varepsilon\to 0$ operators $V^{N,T,S,t}_{t-}$ strongly
converge to a limit $\hat V_{t-}$. $F\hat V_{t-}$ is the operator of
multiplication by the function
$$\sqrt
     \frac{(\tilde T-\tilde t-\tilde S_0+\tilde x)(\tilde t+\tilde N-\tilde x)}{(\tilde t+\tilde N)(\tilde T+\tilde N-\tilde t)}
            +z\cdot\sqrt\frac{(\tilde S_0+\tilde N-\tilde x)\tilde x}{(\tilde t+\tilde N)(\tilde T+\tilde N-\tilde
            t)}.$$
\end{proposition}
\begin{proof}
 Recall that $c_{t-}^{S,t}(k)=c_{t+}^{S,t-1}(k)$ and
 $c_{t-}^{S,t}(k)$ is a real number. Consequently,
 $V^{N,T,S,t}_{t-}=(V^{N,T,S,t-1}_{t+})^*$. Therefore, operators $V^{N,T,S,t}_{t-}$ tend to a limit
$\hat V_{t-}$ and $F\hat V_{t-}=(F\hat V_{t+})^*$. Finally, observe
that if $A$ is the operator of multiplication by $z^{-1}$ acting in
$L_2(S^1)$, then $A^*$ is the operation of multiplication by $z$.
 Thus, $F\hat
V_{t-}$ is given by the desired formula.
\end{proof}

\begin{proposition}
\label{Proposition_V_S+_conv} When $\varepsilon\to 0$ operators
$V^{N,T,S,t}_{S+}$ strongly converge to a limit $\hat V_{S+}$.
$F\hat V_{S+}$ is the operator of multiplication by the function
$$\sqrt
     \frac{(\tilde T-\tilde t-\tilde S_0+\tilde x)(\tilde S_0+\tilde N-\tilde x)}{(\tilde S_0+\tilde N)(\tilde T+\tilde N-\tilde S_0)}
            +z^{-1}\cdot\sqrt\frac{(\tilde t+\tilde N-\tilde x)\tilde x}{(\tilde S_0+\tilde N)(\tilde T+\tilde N-\tilde
            S_0)},$$
\end{proposition}
\begin{proof}
Results of Section \ref{Section_Four_families} imply that
$V^{N,T,S,t}_{S+}\leftrightarrow V^{N,T,t,S}_{t+}$  under the
$S\leftrightarrow t$. Perform this involution, then send
$\varepsilon$ to zero and then switch $S\leftrightarrow t$ again.
Proposition \ref{Proposition_V_t+_conv} implies the result.
\end{proof}

\begin{proposition}
When $\varepsilon\to 0$ operators $V^{N,T,S,t}_{S-}$ strongly
converge to a limit $\hat V_{S-}$. $F\hat V_{S-}$ is the operator of
multiplication by the function
$$\sqrt
     \frac{(\tilde T-\tilde t-\tilde S_0+\tilde x)(\tilde S_0+\tilde N-\tilde x)}{(\tilde S_0+\tilde N)(\tilde T+\tilde N-\tilde S_0)}
            +z\cdot\sqrt\frac{(\tilde t+\tilde N-\tilde x)\tilde x}{(\tilde S_0+\tilde N)(\tilde T+\tilde N-\tilde
            S_0)},$$
\end{proposition}
\begin{proof}
 Same argument as in Proposition \ref{Proposition_V_S+_conv}.
\end{proof}

Now we proceed to the proof of Theorem \ref{Th_bulk_limit}. The
correlation kernel $K(r_i,t_i,x_i;r_j,t_j,x_j)$ defines a $\mathbb
Z\times\mathbb Z$ matrix or, equivalently, an operator acting in
$l_2(\mathbb Z)$ by
$$
 {\cal K}(x,y)=\begin{cases} K(r_i,t_i,x;r_j,t_j,y)
        &x\in{\mathfrak X_{N,T}^{S(r_i),t_i}},y\in{\mathfrak X_{N,T}^{S(r_j),t_j}},\\
            0 &\text{for other } x,y.
  \end{cases}
$$

First, suppose that $r_i<r_j$ or $r_i=r_j$, $t_i>t_j$. In this case
${\cal K}(x,y)$ can be decomposed in the following way:
\begin{multline}
\label{Decomposition} {\cal K}(x,y)={\mathbf
P'}^{N,T,S(r_j),t_j}\cdot V^{N,T,S(r_j),t_j+1}_{t-}\cdot
V^{N,T,S(r_j),t_j+2}_{t-}\dots
V^{N,T,S(r_j),t_i}_{t-}\\
\times V^{N,T,S(r_j-1),t_i}_{S\epsilon_{r_j}}\cdot
V^{N,T,S(r_j-2),t_i}_{S\epsilon_{r_j-1}}\dots
V^{N,T,S(r_i),t_i}_{S\epsilon_{r_i+1}}
\end{multline}
This relation readily follows from the definition of ${\cal K}(x,y)$
and orthogonality relations on functions $\Psi_{l}^{S,t}(x)$.

Observe that norms of all factors in \eqref{Decomposition} are
bounded by $1$. ( ${\mathbf P'}^{N,T,S(r_2),t_2}$ are ortoprojection
operators; norms of $V$'s are bounded since constants $c$ used in
their definition are bounded by $1$). Consequently, convergence of
each factor in \eqref{Decomposition} as $\varepsilon\to 0$ implies
strong convergence of ${\cal K}(x,y)$ to the limit operator
$$
{\cal \hat K}={\mathbf {\hat P'}}(\hat V_{t-})^{t_i-t_j}\cdot {\hat
V}_{S\epsilon_{r_j}}\cdot {\hat V}_{S\epsilon_{r_j-1}}\dots {\hat
V}_{S\epsilon_{r_i+1}} .
$$
(Indeed, multiplication of operators is jointly continuous on
bounded sets in strong operator topology.) The image of ${\cal \hat
K}$ under the Fourier transform is given by
$$
F{\cal \hat K}=F{\mathbf {\hat P'}}(F\hat
V_{t-})^{t_i-t_j}\cdot\prod\limits_{k=r_i}^{r_j-1}F\hat
V_{S\epsilon_{k+1}}.
$$

Performing the inverse Fourier transform we obtain the formula
\begin{multline*}
\lim\limits_{\varepsilon\to 0}
K(r_i,t_i(\varepsilon),x_i(\varepsilon);r_j,t_j(\varepsilon),x_j(\varepsilon))={\cal
\hat K}(\varkappa_i,\varkappa_j)\\= \frac{1}{2\pi
i}\int_{e^{-i\phi}}^{e^{i\phi}}
 \left(\sqrt
     \frac{(\tilde T-\tilde t-\tilde S_0+\tilde x)(\tilde t+\tilde N-\tilde x)}{(\tilde t+\tilde N)(\tilde T+\tilde N-\tilde t)}
            +z\cdot\sqrt\frac{(\tilde S_0+\tilde N-\tilde x)\tilde x}{(\tilde t+\tilde N)(\tilde T+\tilde N-\tilde
            t)}\right)^{t_i-t_j}\\
\times\prod\limits_{k=r_i+1}^{r_j} \left(\sqrt
     \frac{(\tilde T-\tilde t-\tilde S_0+\tilde x)(\tilde S_0+\tilde N-\tilde x)}{(\tilde S_0+\tilde N)(\tilde T+\tilde N-\tilde S_0)}
            +z^{-\epsilon_k}\sqrt\frac{(\tilde t+\tilde N-\tilde x)\tilde x}{(\tilde S_0+\tilde N)(\tilde T+\tilde N-\tilde
            S_0)}\right)\\
  \times\frac{dz}{z^{x_i-x_j+1}}=const_1^{t_i-t_j}\cdot const_2^{r_j-r_i}\cdot K^\bu_{ij}
\end{multline*}
with
$$
 const_1=\sqrt
     \frac{(\tilde T-\tilde t-\tilde S_0+\tilde x)(\tilde t+\tilde N-\tilde x)}{(\tilde t+\tilde N)(\tilde T+\tilde N-\tilde
     t)},
$$
$$
 const_2=\sqrt
     \frac{(\tilde T-\tilde t-\tilde S_0+\tilde x)(\tilde S_0+\tilde N-\tilde x)}{(\tilde S_0+\tilde N)(\tilde T+\tilde N-\tilde
     S_0)}.
$$
 (The integration is performed over the left side of the unit
circle.)

The case $r_i\ge r_j$, $t_i\le t_j$ is quite similar. The only
difference is that instead of the operators $V^{N,T,S,t}_{t-}$ and
$V^{N,T,S,t}_{S\pm}$ we have to use in some sense inverse operators.
Actually, these operators are not invertible and limit operators
$\hat V_{S\pm}$ and $\hat V_{t-}$ might be non-invertible too, but
the difficulties can be avoided if we restrict all operators on the
images of ${\mathbf P}^{N,T,S,t}$ and $\mathbf {\hat P}$. Details of
this trick can be found in Section 3.3 of \cite{Gor}.

The answer for $r_i\ge r_j$, $t_i\le t_j$ is
\begin{multline*}
\lim\limits_{\varepsilon\to 0}
K(r_i,t_i(\varepsilon),x_i(\varepsilon);r_j,t_j(\varepsilon),x_j(\varepsilon))\\=
\frac{1}{2\pi i}\int_{e^{-i\phi}}^{e^{i\phi}}
 \left(\sqrt
     \frac{(\tilde T-\tilde t-\tilde S_0+\tilde x)(\tilde t+\tilde N-\tilde x)}{(\tilde t+\tilde N)(\tilde T+\tilde N-\tilde t)}
            +z\cdot\sqrt\frac{(\tilde S_0+\tilde N-\tilde x)\tilde x}{(\tilde t+\tilde N)(\tilde T+\tilde N-\tilde
            t)}\right)^{t_i-t_j}\\
\times\prod\limits_{k=r_j+1}^{r_i} \left(\sqrt
     \frac{(\tilde T-\tilde t-\tilde S_0+\tilde x)(\tilde S_0+\tilde N-\tilde x)}{(\tilde S_0+\tilde N)(\tilde T+\tilde N-\tilde S_0)}
            +z^{-\epsilon_k}\sqrt\frac{(\tilde t+\tilde N-\tilde x)\tilde x}{(\tilde S_0+\tilde N)(\tilde T+\tilde N-\tilde
            S_0)}\right)^{-1}\\
  \times\frac{dz}{z^{x_i-x_j+1}}=const_1^{t_i-t_j}\cdot
  const_2^{r_j-r_i}\cdot
  K^\bu_{ij},
\end{multline*}
where $const_1$ and $const_2$ are as above. (The integration is
performed over the right side of the unit circle.)

Since in the determinant for correlation functions the prefactors
$const_1^{t_i-t_j}const_2^{r_j-r_i}$ cancel out, the proof is
complete.

\end{proof}

\end{document}